\documentstyle[12pt,leqno,twoside]{article}
\input{amssym.def}
\title{ \LARGE \rm Almost Quaternion-Hermitian Manifolds}
\author{\normalsize {\rm FRANCISCO MART{\'I}N  CABRERA}}
\date{}
\newcounter{sec}
\newtheorem{proposition}{\sc PROPOSITION}[sec]
\newtheorem{lemma}[proposition]{\sc LEMMA}

\newtheorem{observation}[proposition]{\sc REMARK}
\newcommand{\Dt}{{\it Proof}.- }

\newcommand{\com}{\makebox[7pt]{\raisebox{1.5pt}{\tiny$\circ$}}}

\newcommand{\sumcic}{\mathop{\mbox{\large {$\frak S$}\vrule width 0pt
depth 2pt}}}

\newcommand{\seccion}[1]{\setcounter{equation}{0}\addtocounter{sec}{1}
   \section{\normalsize \bf \arabic{sec}. #1  \hspace{3em} }  }

\def\thebibliography#1{\section*{\normalsize\bf
References}
        \list {\rm \arabic{enumi}.}{\settowidth\labelwidth{[#1]}
                \leftmargin\labelwidth  \advance\leftmargin\labelsep
                 \usecounter{enumi}}
         \def\newblock{\hskip .11em plus .33em minus .07em}
         \sloppy\clubpenalty4000\widowpenalty4000
         \sfcode`\.=1000\relax}
\newcommand{\cqd}{\hfill $\Box$\par}
\newcommand{\qedhere}{\mbox{\qquad$\Box$\hss}}

\begin{document}

\baselineskip=5mm
\pagestyle{myheadings}

\maketitle \vspace{-10mm}

\noindent {\footnotesize {\rm Departmento de Matem{\'a}tica
Fundamental, Universidad
 de La Laguna, 38200 La Laguna, Tenerife, Spain. e-mail: {\tt
 fmartin@ull.es}}}

\begin{abstract}
{\footnotesize \noindent {\bf Abstract.} Following the point of
view of Gray and Hervella, we derive detailed  conditions which
characterize each one of the classes of almost
quaternion-Hermitian $4n$-manifolds, $n>1$. Previously, by
completing a basic result of A. Swann, we give explicit
descriptions of the tensors contained in  the space of covariant
derivatives of the fundamental form $\Omega$ and split the
coderivative of $\Omega$ into its ${\sl Sp}(n){\sl
Sp}(1)$-components. For $4n>8$, A. Swann also proved that all the
information about  the intrinsic torsion $\nabla \Omega$ is
contained in the exterior derivative ${\sl d} \Omega$. Thus, we
give alternative conditions, expressed in terms of ${\sl d}
\Omega$, to characterize the different classes of almost
quaternion-Hermitian manifolds.}
\end{abstract}




\baselineskip=5mm


\seccion{Introduction} \setcounter{proposition}{0} \noindent An
almost quaternion-Hermitian manifold is a Riemannian $4n$-manifold
($n >1$) which admits a ${\sl Sp}(n){\sl Sp}(1)$-structure, i.e.,
 a reduction of
its frame bundle to the subgroup ${\sl Sp}(n){\sl
Sp}(1)$ of ${\sl SO}(4n)$.
 These  manifolds are of special interest because
 ${\sl
Sp}(n){\sl Sp}(1)$ is included  in the list of Berger
(\cite{B})  of
possible holonomy groups
of locally irreducible Riemannian manifolds that are not locally
symmetric. An almost quaternion-Hermitian manifold is
said to be quaternion-K{\"a}hler, if its reduced holonomy
group is a subgroup of ${\sl Sp}(n){\sl Sp}(1)$.

A general approach to classify $G$-structures compatible with the
metric in a Riemannian manifold is described as follows. Given any
Riemannian manifold $M^m$ and a Lie group $G$ that is the
stabilizer of some tensor $\varphi$ on ${\Bbb R}^m$, that is, $ G
= \{ g \in {\sl SO}(m)\, : \, g \varphi =  \varphi \}$, a
$G$-structure on $M$ defines a global tensor $\varphi$ on $M$, and
it can be shown that $\nabla \varphi$ (the so-called intrinsic
torsion of the $G$-structure, where $\nabla$ denotes the
Levi-Civita connection) is a section of the vector bundle ${\cal
W} = {\sl T}^{\ast} M \otimes {\frak g}^{\perp}$, being ${\frak
s}{\frak o}(m) = {\frak g} \oplus {\frak g}^{\perp}$. The action
of $G$ splits ${\cal W}$ into irreducible components, say ${\cal
W} = {\cal W}_1 \oplus \ldots \oplus {\cal W}_k$. Then,
$G$-structures on $M$ can be classified in at most $2^k$ classes.

This systematic way of classifying $G$-structures was initiated by
A. Gray and L. Hervella in \cite{GH}, where they considered the
case $G = {\sl U}(n)$ (almost Hermitian structures), giving rise
to sixteen classes of almost Hermitian manifolds.  Later, diverse
authors have studied the situation for other $G$-structures: ${\sl
G}_2$, ${\sl Spin}(7)$, ${\sl Spin}(9)$ etc.
 In the
present paper we are concerned with the case $G = {\sl Sp}(n){\sl
Sp}(1)$, i.e., we consider a $4n$-manifold  $M$ with an ${\sl
Sp}(n){\sl Sp}(1)$-structure, $n>1$. Associated with this type of
geometric structure there is a four-form $\Omega$, given by
\ref{quafor},  such that $\Omega^n$ is nowhere zero. An
orientation on the manifold can be defined regarding a constant
multiple of $\Omega^n$ as volume form, and ${\sl Sp}(n){\sl
Sp}(1)$ is the stabilizer of $\Omega$ under the action of ${\sl
SO}(4n)$. Moreover, $\nabla \Omega$ is contained in ${\cal W} =
{\sl T}^{\ast} M \otimes \left( {\frak s}{\frak p}(n)\oplus {\frak
s}  {\frak p}(1) \right)^{\perp}$. The splitting of the ${\sl
Sp}(n){\sl Sp}(1)$-module ${\cal W}$ was shown by A. Swann in
\cite{Sw1} (see Proposition \ref{swann1} of the present paper).
Such a decomposition consists of four or six irreducible
components, for $4n=8$ or $4n
> 8 $, respectively. Therefore, it is already known how many
classes of almost quaternion-Hermitian manifolds we have at most,
$16$ or $64$. But actually, to our knowledge, except some partial
results (see \cite{Sw4}), there are not explicit descriptions of
each one of such classes. Thus, the main purpose of the present
paper is to show conditions which define each class, and that is
the content of the second column of Table 2 (page
\pageref{tados}).

For some $G$-structures, the intrinsic torsion $\nabla \varphi$ is
totally determined by tools of the exterior algebra. This implies
that we would have an alternative easier way to characterize  each
class. For an ${\sl Sp}(n){\sl Sp}(1)$-structure, when $4n
> 8 $, this is such a case. A.  Swann in \cite{Sw1} proved that
all the information about the covariant derivative $\nabla \Omega$
is contained in the exterior derivative ${\sl d} \Omega$. Thus
alternative conditions expressed in terms of ${\sl d} \Omega$ are
displayed in the third column of Table 2 to characterize the
different classes of almost quaternion-Hermitian manifolds of
dimension higher than eight. Since for dimension eight, ${\sl
d}\Omega$ only contains partial information of $\nabla \Omega$, in
this case  we can only characterize, via ${\sl d} \Omega$, eight
of the possible sixteen classes. Such characterizations are shown
in Table 3.


\seccion{Almost quaternion-Hermitian structures}
\setcounter{proposition}{0} \baselineskip=5mm \noindent Almost
quaternion-Hermitian manifolds have been broadly treated by
diverse authors (\cite{S1},\cite{Sw1}, etc.). In this section we
recall basic definitions and known results.

     A $4n$-dimensional manifold $M$ $(n > 1)$ is said to be {\it
almost quaternion-Hermitian}, if  $M$ is equipped with a
Riemannian metric $\langle , \rangle$ and a rank-three subbundle
${\cal G}$ of the endomorphism bundle
 ${\sl End} \; TM$,  such that locally ${\cal G}$  has an {\it adapted
basis} $I, J, K$ with $I^2=J^2 = -1$, $K=IJ = -JI$, and $\langle A
X,A Y \rangle = \langle X,Y \rangle $, for $A =I,J,K$. This is
equivalent to saying that $M$ has a reduction of its structure
group to ${\sl Sp}(n){\sl Sp}(1)$.


There are the three local K{\"a}hler-forms $\omega_{A} (X,Y) = \langle
X ,A Y \rangle$, $A=I,J,K$. From these one may define a global
four-form, non degenerate $\Omega$, the {\it fundamental form}, by
the local formula
\begin{equation}
\Omega = \omega_I \wedge \omega_I + \omega_J \wedge \omega_J +
\omega_K \wedge \omega_K. \label{quafor}
\end{equation}

Now, we recall some facts about almost quaternion-Hermitian
manifolds in relation with representation theory. We will  follow
the ${\sl E}$-${\sl H}$-formalism  used in \cite{S1},
\cite{Sw1}-\cite{Sw4} and we refer to \cite{BtD} for general
information on representation theory. Thus, ${\sl E}$ is the
fundamental representation of ${\sl Sp}(n)$ on ${\Bbb C}^{2n}
\cong {\Bbb H}^n$ via left multiplication by quaternionic
matrices, considered in $\mbox{\rm GL}(2n, {\Bbb C})$, and ${\sl
H}$ is the representation of ${\sl Sp}(1)$ on ${\Bbb C}^2 \cong
{\Bbb H}$ given by $q . \zeta = \zeta \overline{q}$, for $q \in
{\sl Sp}(1)$ and $\zeta \in {\sl H}$. An ${\sl Sp}(n){\sl
Sp}(1)$-structure on a manifold $M$ gives rise to local bundles
${\sl E}$ and ${\sl H}$ associated to these representation and
identifies ${\sl T} M \otimes_{\Bbb R} {\Bbb C} \cong {\sl E}
\otimes_{\Bbb C} {\sl H}$.

On ${\sl E}$, there is an ${\sl Sp}(n)$-invariant complex
symplectic form $\omega_{\sl E}$ and a Hermitian inner product
given by $\langle x , y \rangle_{\Bbb C} = \omega_{\sl E} (x ,
\widetilde{y})$, for all $x,y \in {\sl E}$ and being
$\widetilde{y} = j y $ ($y \to \widetilde{y}$ is a quaternionic
structure map on ${\sl E}={\Bbb C}^{2n}$ considered as left
complex vector space). The mapping $x \to x^{\omega} = \omega_{\sl
E} ( \cdot, x)$ gives us an identification of ${\sl E}$ with its
dual ${\sl E}^{\ast}$. If $\{ e_1, \dots , e_n, \widetilde{e}_1,
\dots , \widetilde{e}_n \}$ is a complex orthonormal basis for
${\sl E}$, then $ \omega_{\sl E} = e^{\omega}_i \wedge
\widetilde{e}^{\omega}_i = e^{\omega}_i \widetilde{e}^{\omega}_i -
\widetilde{e}^{\omega}_i e^{\omega}_i, $ where we have used the
summation convention and omitted tensor product signs. These
conventions will be used throughout the paper.

The ${\sl Sp}(1)$-module ${\sl H}$ will be also considered as left
complex vector space. Regarding ${\sl H}$ as  4-dimensional real
space with the Euclidean metric $\langle , \rangle$ such that $\{
1 , i , j , k \}$ is an orthonormal basis, the complex symplectic
form $\omega_{\sl H}$  is given  by $\omega_{\sl H} = 1^{\#}
\wedge j^{\#} + k^{\#} \wedge i^{\#} + i (1^{\#} \wedge k^{\#} +
i^{\#} \wedge j^{\#})$, where $h^{\#}$ is the real one-form given
by $q \to \langle h , q \rangle$. We also have the identification,
$q \to q^{\omega} = \omega_{\sl H} (\cdot , q )$, of $H$ with its
dual ${\sl H}^{\ast}$  as complex space.  On ${\sl H}$, we have a
quaternionic structure map given by $q = z_1 + z_2 j \to
\widetilde{q} = jq = - \overline{z}_2 + \overline{z}_1 j$, where
$z_1,z_2 \in {\Bbb C}$ and $\overline{z}_1, \overline{z}_2 $ are
their conjugates. If $h \in {\sl H}$ is such that $\langle h , h
\rangle =1$, then $\{ h , \widetilde{h} \}$ is a basis of the
complex vector space ${\sl H}$ and $\omega_{\sl H} = h^{\omega}
\wedge \widetilde{h}^{\omega}$.

The irreducible representations of ${\sl Sp}(1)$ are the symmetric
powers ${\sl S}^k {\sl H} \cong {\Bbb C}^{k+1}$. An irreducible
representation of ${\sl Sp}(n)$ is determined by its dominant
weight $(\lambda_1,\dots,$ $ \lambda_n)$, where $\lambda_i$ are
integers with $\lambda_1 \geq \lambda_2 \geq \dots \geq \lambda_n
\geq 0$. This representation will be denoted by $V^{ (\lambda_1,
\dots , \lambda_r)}$, where $r$ is the largest integer such that
$\lambda_r > 0$. Familiar notation is used for some of these
modules, when possible. For instance,
 $V^{(k)} = {\sl S}^k {\sl E}$,  the $k$th symmetric power of ${\sl E}$, and  $V^{(1,
\dots , 1)} = \Lambda_0^r {\sl E}$, where there are $r$ ones in
exponent and $\Lambda_0^r {\sl E}$ is the ${\sl Sp}(n)$-invariant
complement to $\omega_E \Lambda^{r-2} {\sl E}$ in $\Lambda^r {\sl
E}$. Also ${\sl K}$ will denote the module $V^{(21)}$, which
arises in the decomposition $ {\sl E} \otimes \Lambda_0^2 {\sl E}
\cong \Lambda_0^3 {\sl E} + {\sl K} + {\sl E}$, where $+$ denotes
direct sum.

As we have pointed out in the introduction an analogue of the
Gray-Hervella classification may be obtained for almost quaternion
hermitian manifolds by considering the covariant derivative
$\nabla \Omega$ of the fundamental $4$-form $\Omega$
(\ref{quafor}). For dimension at least $12$, this leads to $2^6
=64$ classes. Such quaternionic classes come from the following
${\sl Sp}(n){\sl Sp}(1)$-decomposition.
\begin{proposition}[Swann~\cite{Sw1}] \label{swann1}
{\it  The covariant derivative of the fundamental form $\Omega$ of
an almost
  quaternion-Hermitian manifold $M$ of dimension at least $8$, has the
  property
  \begin{displaymath}
  \nabla \Omega \in T^* M \otimes \left( {\frak s}{\frak
p}(n)\oplus {\frak s}  {\frak p}(1) \right)^{\perp}  \cong T^* M
\otimes \Lambda_0^2{\sl E} {\sl S}^2 {\sl H} = ( \Lambda_0^3 {\sl
E} +
  {\sl K} + {\sl E}) ( S^3 {\sl H} + {\sl H}).\qedhere
  \end{displaymath}}
\end{proposition}

\noindent If the dimension of $M$ is at least~$12$, all the
modules of the sum are non-zero.  For an eight-dimensional
manifold $M$, we have $\Lambda_0^3 E S^3 H = \Lambda_0^3 E H = \{
0 \}$.
\begin{observation}
{\it Let $V$ be a complex $G$-module equipped with a real
structure, where $G$ is a Lie group. In the previous proposition
and most of the times throughout this paper,
 $V$ also denotes the  $(+1)$-eigenspace of the structure map which
 is a real $G$-module. The
context should  tell us which space are referring to. However, if
there is   risk of confusion or when we feel that a clearer
exposition is needed, we denote the second mentioned space by
$[V]$. }
\end{observation}

For sake of simplicity, we are referring to a real vector  space
${\cal V}$  of dimension $4n$ with $n>1$. In the final
conclusions, we will apply the obtained results to an almost
quaternion-Hermitian manifold $M$ and   ${\cal V}$ will be $T_p
M$, the tangent space at $p \in M$.

Thus, ${\cal V}$ is  equipped with an Euclidean metric $\langle ,
\rangle$ and a three-dimensional subspace  ${\cal G}$ of the
endomorphism space
 ${\sl End} \; {\cal V}$ such that  ${\cal G}$  admits a
basis $I,J,K$, with $I^2 = J^2 = -1$, $K=IJ= - JI$ and $\langle A
x,A y \rangle = \langle x,y \rangle $, for $A=I,J,K$. From  the
two-forms $\omega_{A} (x,y) = \langle x ,A y \rangle$, a four-form
$\Omega$ is defined as in \ref{quafor}.

Moreover, we can consider ${\cal V} $ as a left complex vector
space by defining $ (a+bi) x =  a x + b I x$, $x \in {\cal V}$ and
$a+bi \in {\Bbb C}$.  We denote  ${\sl E}$ when we refer to ${\cal
V} $ as complex vector space and we have ${\sl dim} E = {\sl
dim}_{\Bbb C} {\cal V}= 2n$. The ${\sl Sp}(n)$-representation
${\sl E}$ has a quaternionic structure map ${\sl E} \to {\sl E}$
given by $x \to \widetilde{x} =  J x$ and is equipped with the
${\sl Sp}(n)$-invariant complex symplectic form $\omega_{\sl E} =
- \omega_J - i \omega_K$.  Note that $x^{\omega} = - J x - i K x$,
where $x$ also denote the dual form $y \to \langle y , x \rangle$,
for all $y \in {\cal V}$(the identification ${\cal V} = {\cal
V}^{\ast}$).

The complexification of the real vector space ${\cal V}$ can be
identified with  ${\sl E} \otimes_{\Bbb C} {\sl H}$ via the
isomorphism $ {\cal V} \otimes_{\Bbb R} {\Bbb C}  \to   {\sl E}
\otimes_{\Bbb C} {\sl H}$, defined by $
  x \otimes_{\Bbb R} z  \to  x \otimes_{\Bbb C} z h + J x
\otimes_{\Bbb C} z \widetilde{h}$, where we have fixed $h \in {\sl
H}$ such that $\langle h , h \rangle=1$. On ${\sl E} \otimes_{\Bbb
C} {\sl H}$, we can consider the real structure map $x
\otimes_{\Bbb C} q \to \widetilde{x} \otimes_{\Bbb C}
\widetilde{q}$, which corresponds with the map $x \otimes_{\Bbb R}
z \to  x \otimes_{\Bbb R} \overline{z}$ on ${\cal V} \otimes_{\Bbb
R} {\Bbb C}$. The $(+1)$-eigenspace $[{\sl EH}] \cong {\cal V}
\otimes_{\Bbb R} {\Bbb R} \cong {\cal V}$ is the real irreducible
${\sl Sp}(n){\sl Sp}(1)$-representation determined by the complex
irreducible ${\sl Sp}(n){\sl Sp}(1)$-representation ${\sl E}
\otimes_{\Bbb C} {\sl H}$ of real type.

Taking $\langle,\rangle_{\Bbb C} = \langle,\rangle + i \omega_I$
into account, for $x \in {\sl E}$, we obtain that the restrictions
of $ x^{\omega} h^{\omega}$, $ x^{\omega} \widetilde{h}^{\omega}$,
$\widetilde{x}^{\omega} h^{\omega}$   and $ \widetilde{x}^{\omega}
\widetilde{h}^{\omega}$ to $ {\cal V} \cong [{\sl EH}]$ are given
by
\begin{equation}
\begin{array}{lll} \label{equiache}
 x^{\omega} h^{\omega}_{|{\cal V}}  =  x - i I x, &
 \quad &
 x^{\omega} \widetilde{h}^{\omega}_{|{\cal V}}  =  - Jx - i
K x, \\
\widetilde{x}^{\omega} h^{\omega}_{|{\cal V}}  =  J x - i K x, &
\quad & \widetilde{x}^{\omega} \widetilde{h}^{\omega}_{|{\cal V}}
= x + i I x.
\end{array}
\end{equation}

\setlength{\baselineskip}{5mm} The following conventions will be
used in the sequel.  If $b$ is a $(0,s)$-tensor, we write
\begin{displaymath}
  \begin{array}{c}
    A_{(i)}b(X_1, \dots, X_i, \dots , X_s) = - b(X_1, \dots , AX_i, \dots ,
    X_s),\\
    A b(X_1,\dots,X_s) = (-1)^sb(AX_1,\dots,AX_s),\\
{\sl i}_{A}b   =  (A_{(1)} + \dots + A_{(s)})b,
  \end{array}
\end{displaymath}
for $A=I,J,K$.  We also consider the extension of~$\langle\cdot
,\cdot \rangle$ to $(0,s)$-tensors given by
\begin{displaymath}
  \langle a, b \rangle = \mbox{$\frac{1}{s!}$} \, a (
  e_{i_1}, \dots , e_{i_s})  b (e_{i_1}, \dots , e_{i_s}),
\end{displaymath}
where $\{ e_1, \dots , e_{4n} \}$ an orthonormal basis for ${\cal
V}$. Finally, we define the ${\sl Sp}(n){\sl Sp}(1)$-map $L\, : \,
\Lambda^p {\cal V}^* \to \Lambda^p {\cal V}^*$  by
\begin{equation} \label{lmant}
L(b) =  \sum_{A=I,J,K} \sum_{1\leq i < j \leq p} A_{(i)}A_{(j)}b.
\end{equation}

\seccion{The space of covariant derivatives of $\Omega$}{\indent}
\setcounter{proposition}{0} In this section we will give an
explicit description of the tensors which are contained in the
space of covariant derivatives of the fundamental  four form
$\Omega$. For such a purpose, we consider the space $ {\cal
V}^{\ast} \otimes \left( {\frak s}{\frak p}(n) + {\frak s} {\frak
p}(1) \right)^{\perp} \cong {\cal V}^{\ast} \otimes \Lambda_0^2E
{\sl S}^2 H \subseteq {\cal V}^{\ast} \otimes \Lambda^2 {\cal
V}^{\ast}$ consisting  of the tensors $c$ which, for an adapted
basis $I,J,K$
 of ${\cal G}$, satisfy:
\begin{enumerate}
\item[i)$^*$] $c + I_{(2)}I_{(3)}  c + J_{(2)} J_{(3)}  c +
K_{(2)}K_{(3)} c = 0$;
 \item[ii)$^*$] $\langle \cdot \lrcorner c , \omega_{A}
\rangle = 0$, for $A = I,J,K$.
\end{enumerate}
We recall that $\Lambda^2 {\cal V}^{\ast} = S^2 {\sl E} + {\sl
S}^2 {\sl H} +  \Lambda_{0}^2 {\sl E} {\sl S}^2 {\sl H} $, where
$S^2 {\sl E} \cong {\frak s}{\frak p}(n)$ and ${\sl S}^2 {\sl H}
\cong {\frak s}{\frak p}(1)$ are the Lie algebras of ${\sl Sp}(n)$
and ${\sl Sp}(1)$, respectively.

Now  we define the
  ${\sl Sp}(n){\sl Sp}(1)$-map
\begin{displaymath}
\begin{array}{rcl}
 {\cal F} :  {\cal V}^{\ast} \otimes \Lambda_0^2{\sl E} {\sl S}^2 {\sl H}  &
 \rightarrow &
{\cal V}^{\ast} \otimes \Lambda^4 {\cal V}^{\ast}, \\
 c & \rightarrow & \mbox{$\frac{1}{4}$} \sum_{A=I,J,K} {\sl i}_{A} (\cdot \lrcorner c) \wedge
\omega_{A}.
\end{array}
\end{displaymath}
Since $\Lambda_0^2 {\sl E} {\sl S}^2 {\sl H}$ is an irreducible
${\sl Sp}(n){\sl Sp}(1)$-module, by    Schur's Lemma (see \cite[p.
110]{BtD}),
 $
 {\cal F} : {\cal V}^{\ast} \otimes
\Lambda_0^2{\sl E} {\sl S}^2 {\sl H} \to \, {\cal F}({\cal
V}^{\ast} \otimes \Lambda_0^2{\sl E} {\sl S}^2 {\sl H}) $ is an
${\sl Sp}(n){\sl Sp}(1)$-isomorphism. If $a= {\cal F}(c)$, one can
check
\begin{equation} \label{ace}
-8 n c (x, y , z ) = \langle x \lrcorner a , y \wedge (z \lrcorner
\Omega) - z \wedge (y \lrcorner \Omega) \rangle,
\end{equation}
for all $x, y,z \in {\cal V}$. Thus the map $a \to c$, defined by
\ref{ace}, is the inverse map of ${\cal F}$.

\begin{proposition} \label{Perpsp}
{\it If  $a \in {\cal V}^{\ast} \otimes \Lambda^4 {\cal
V}^{\ast}$, the following conditions are equivalent:
\begin{enumerate}
\item[1)] $a \in {\cal V}^{\ast} \otimes \Lambda_0^2{\sl E} {\sl
S}^2 {\sl H}$.
 \item[2)] There exists an adapted basis $I, J, K$ of ${\cal
G}$ and  a unique triplet $c_I, c_J, c_K \in {\cal V}^{\ast}
\otimes \Lambda^2 {\cal V}^{\ast}$ associated with such a basis,
satisfying $a = \sum_{A=I,J,K} c_{A} \wedge \omega_{A}$ and the
conditions:
\begin{enumerate}
\item[{\rm i)}] $A_{(2)}A_{(3)} c_{A} = - c_{A}$, for $A=I,J,K$.
\item[{\rm ii)}] $I_{(2)}J_{(3)} c_{K} + J_{(2)}K_{(3)} c_{I} +
K_{(2)}I_{(3)} c_{J} =0$.
 \item[{\rm iii)}] $\langle \cdot \lrcorner c_{A} , \omega_{B} \rangle = 0$, for $A,B=I,J,K$.
\end{enumerate}
\item[3)] There exists an adapted basis $I,J,K$ of ${\cal G}$ and
three tensors $d_I, d_J, d_K \in {\cal V}^{\ast} \otimes \Lambda^2
{\cal V}^{\ast}$, satisfying $ a = \sum_{A=I,J,K} d_{A} \wedge
\omega_{A}$ and the conditions:
\begin{enumerate}
\item[{\rm i)$'$}] $A_{(2)}A_{(3)} d_{A} = - d_{A}$, for
$A=I,J,K$.
 \item[{\rm ii)$'$}] $I_{(2)}J_{(3)} d_{K} +
J_{(2)}K_{(3)} d_{I} + K_{(2)}I_{(3)} d_{J} =0$.
\end{enumerate}
\end{enumerate} }
\end{proposition}
\Dt If $a \in {\cal V}^{\ast} \otimes \Lambda_0^2E {\sl S}^2 H$,
then $a = {\cal F}(c)$. Therefore, $a = \sum_{A=I,J,K} c_{A}
\wedge \omega_{A}$, where $c_A = \frac{1}{4} {\sl i}_{A} (\cdot
\lrcorner c)$. Taking conditions i)$^*$ and ii)$^*$ into account,
it follows that $c_I$, $c_J$ and $c_K$ satisfy conditions i), ii)
and iii) of the lemma.

On the other hand, if $ a = \sum_{A=I,J,K} c_{A} \wedge
\omega_{A}$ with $c_I$, $c_J$, $c_K$ satisfying the conditions i),
ii) and iii), we consider the tensor  $ c = \frac{1}{2}
\sum_{A=I,J,K} {\sl i}_{A} (\cdot \lrcorner c_{A})$. Taking i),
ii) and iii) into account, it follows that $c$ satisfies
conditions i)$^*$ and ii)$^*$ and ${\cal F}(c) =a$. Hence $a \in
{\cal V}^{\ast} \otimes  \Lambda_0^2 E {\sl S}^2H$. Hence, 1) is
equivalent to 2). The triplet  $c_I$, $c_J$ and $c_K$ is unique,
because one can check  $- 4n c_A(x,y,z) = a(x,y,z,e_r, Ae_r)$.

In order to prove 3) implies 2), let us consider  the tensors
$c_{I}$, $c_J$, $c_K$ given by
$$
c_{I} = d_{I} + \mbox{$\frac{1}{2n}$} \langle d_{I} , \omega_{J}
\rangle \omega_{J} + \mbox{$\frac{1}{2n}$} \langle d_{I} ,
\omega_{K} \rangle \omega_{K}.
$$
The tensors $c_J$ and $c_K$ are also defined from this identity by
cyclically permuting $I,J,K$. Throughout this paper, when we have
three equations depending cyclically of $I,J,K$ as here,  we will
write only one equation without further comment.
 Note that $d_{A} \in {\cal V}^{\ast} \otimes ({\sl S}^2 {\sl H} +
\Lambda_{0}^2 {\sl E} {\sl S}^2 {\sl H})$ and $c_{A}$ is its
projection on ${\cal V}^{\ast} \otimes \Lambda_{0}^2 {\sl E} {\sl
S}^2 {\sl H}$. Also note that $ a = \sum_{A=I,J,K} c_{A} \wedge
\omega_{A}$ and that $c_I$, $c_J$ and $c_K$ satisfy i), ii) and
iii). \cqd

\begin{observation}
{\it  In the context of almost quaternion-Hermitian manifolds, we
consider $a= \nabla \Omega$ and we have $d_{A} = 2\nabla
\omega_{A}$. In fact, $2 \nabla \omega_A$ satisfies condition {\rm
i)$'$} (\cite{GH}) and it is an easy exercise  to check that
$\nabla \omega_I$, $\nabla \omega_J$, $\nabla \omega_K$ satisfy
condition {\rm ii)$'$}. Tensors $c_A$ satisfying conditions {\rm
i)}, {\rm ii)} and {\rm iii)} are given by
$$
c_I = 2 \nabla \omega_I + \mbox{$\frac{1}{n}$} \langle \nabla
\omega_I, \omega_J \rangle \omega_J + \mbox{$\frac{1}{n}$} \langle
\nabla \omega_I, \omega_K \rangle \omega_K.
$$
  }
\end{observation}

\begin{lemma} \label{lemma:La2a}
{\it Let ${\sl L} \, : \Lambda^4 {\cal V}^* \to \Lambda^4 {\cal
V}^*$ be  defined by \ref{lmant}, then for all $a \in
\Lambda_0^2{\sl E} {\sl S}^2 {\sl H} \subseteq \Lambda^4 {\cal
V}^* $ we have ${\sl L}(a) = 2a$. In particular, in a context of
almost quaternion-Hermitian manifolds  we have ${\sl L}(\nabla_X
\Omega) = 2 \nabla_X \Omega$, for all tangent vector $X$.
 }
\end{lemma}
\Dt By Proposition \ref{Perpsp}, for all four-form $a \in
\Lambda_0^2{\sl E} {\sl S}^2 {\sl H}$, we have $a = \sum_{A=I,J,K}
c_A \wedge w_A$, where $c_I$, $c_J$, $c_K$ are skew-symmetric
two-forms satisfying conditions derived directly from i), ii) and
iii) of Proposition \ref{Perpsp}. Moreover, $-4n c_A(x,y) =
a(x,y,e_r, Ae_r)$.

By straightforward computation, we obtain
\begin{eqnarray*}
{\sl L}(a)(x,y,e_r,Ie_r) & = & - 4n \left( - 2 c_I (x,y) +
c_I(Jx,Jy) + c_I(Kx,Ky) + 2c_J (Kx,y) \right. \\
&& \left. \qquad \;  + 2c_J(x,Ky) - 2 c_K(Jx,z) - 2 c_K(x,Jy)
\right).
\end{eqnarray*}
Now, making use of the conditions satisfied by $c_I$, $c_J$ and
$c_K$, we will have \linebreak ${\sl L}(a)(x,y,e_r,Ie_r) = - 8n
c_I(x,y)$. In a similar way, we can obtain ${\sl
L}(a)(x,y,e_r,Je_r) = - 8n c_J(x,y)$ and ${\sl L}(a)(x,y,e_r,Ke_r)
= - 8n c_K(x,y)$.

If $a\neq 0$, then $c_I \neq 0$ or $c_J\neq 0$ or $c_K \neq 0$.
Therefore, ${\sl L}(a) \neq 0$. Since there is only one copy of
$\Lambda_0^2{\sl E} {\sl S}^2 {\sl H}$ in the ${\sl Sp}(n){\sl
Sp}(1)$-decomposition of $\Lambda^4 {\cal V}^{\ast}$ (see \cite[p.
127]{S2}), by Schur's Lemma,  ${\sl L}$ is an ${\sl Sp}(n){\sl
Sp}(1)$-isomorphism on $\Lambda_0^2{\sl E} {\sl S}^2 {\sl H}$.
Thus $L(a) = 2 \sum_{A=I,J,K} c_A \wedge \omega_A = 2 a$. \cqd
\vspace{3mm}

There is  an alternative description of  ${\cal V}^{\ast} \otimes
\Lambda_0^2{\sl E} {\sl S}^2 {\sl H} $ as a subspace of ${\cal
V}^{\ast} \otimes \Lambda^4 {\cal V}^{\ast} = {\cal V}^{\ast}
\otimes  \Lambda^4 ({\sl EH}) $. In fact, we have
$
 \Lambda_0^2{\sl E} {\sl S}^2 {\sl H} \cong \left(  \Lambda_0^2{\sl E}  \wedge
\omega_{\sl E} \right) \otimes \Lambda^2   {\sl S}^2 {\sl H}
\subseteq  ({\sl EH})^4$.  Since  the restriction to $\left(
\Lambda_0^2{\sl E} \wedge \omega_{\sl E} \right) \otimes \Lambda^2
{\sl S}^2 {\sl H} $ of the alternation map ${\sl a} : ({\sl
EH})^4$ $ \to$ $ \Lambda^4 ({\sl EH})$ is non-zero, by Schur's
Lemma, we have
\begin{equation} \label{cuaeh}
{\cal V}^{\ast} \otimes \Lambda_0^2{\sl E} {\sl S}^2 {\sl H} \cong
{\cal V}^{\ast} \otimes  {\sl a} \left( \left( \Lambda_0^2{\sl E}
\wedge \{ \omega_{\sl E} \} \right) \otimes \Lambda^2 {\sl S}^2
{\sl H} \right) \subseteq {\cal V}^{\ast} \otimes  \Lambda^4 ({\sl
EH}).
\end{equation}


Next  we obtain an ${\sl Sp}(n){\sl Sp}(1)$-automorphism on ${\cal
V}^{\ast} \otimes \Lambda_0^2{\sl E} {\sl S}^2 {\sl H}$, which
will play significant r{\^o}le in explicit descriptions of classes of
almost quaternion Hermitian manifolds. Thus,  we consider ${\cal
L}
 :   {\cal V}^{\ast} \otimes \Lambda_0^2{\sl E} {\sl S}^2 {\sl
H} \to   {\cal V}^{\ast} \otimes \Lambda_0^2{\sl E} {\sl S}^2 {\sl
H}$, defined by
\begin{equation} \label{eleauto}
{\cal L}(a) =  \sum_{A=I,J,K} A_{(1)}(A_{(2)} + \dots + A_{(5)})
a.
\end{equation}
Relative to the map ${\cal L}$, we have the following results.
\begin{proposition} \label{laele}
{\it
\begin{enumerate}
\item[{\rm i)}] ${\cal L}$ is an ${\sl Sp}(n){\sl
Sp}(1)$-isomorphism.
 \item[{\rm ii)}] $(\Lambda_0^{3} {\sl E} +
{\sl K} + {\sl E}){\sl H}$ consists of    $ a \in {\cal V}^{\ast}
\otimes \Lambda_0^2{\sl E} {\sl S}^2 {\sl H}$  such that ${\cal
L}(a) = 4 a$.
 \item[{\rm iii)}] $(\Lambda_0^{3} {\sl E} + {\sl K} + {\sl E}) {\sl S}^3 {\sl H} $ consists of
 $a \in {\cal V}^{\ast} \otimes \Lambda_0^2{\sl E} {\sl S}^2 {\sl H}$ such that ${\cal L}(a) = -2 a$.
\end{enumerate}
}
\end{proposition}
\Dt Let us consider the  ${\sl Sp}(1)$-maps $Id \otimes
\omega_{23}: {\sl H} \otimes \Lambda^2 {\sl S}^2 {\sl H}  \to {\sl
H} \otimes {\sl S}^2 {\sl H} \cong {\sl S}^3 {\sl H} + {\sl H}$,
$s : {\sl H} \otimes {\sl S}^2 {\sl H}  \to  {\sl S}^3{\sl H}$ and
$\omega_{12} : {\sl H} \otimes {\sl S}^2 {\sl H}   \to {\sl H}$,
where $Id$ is the identity on ${\sl H}$, $s$ is the symmetrization
and $\omega_{\alpha \beta}$ is the $\alpha \beta$- contraction in
$\Lambda^2 ({\sl S}^2 H)$ by $\omega_{\sl H}$. Note that $Id
\otimes \omega_{23}$ is an isomorphism, ${\sl H} \cong {\sl ker}\,
s \com (Id \otimes \omega_{23}) $, ${\sl S}^3 {\sl H} \cong {\sl
ker}\, \omega_{12} \com (Id \otimes \omega_{23}) $ and both are
${\sl Sp}(1)$-subspaces of $ {\sl H} \otimes \Lambda^2 {\sl S}^2
{\sl H}$. Taking this into account, it follows that
$$
\begin{array}{c}
\qquad h^{\omega} \otimes h^{\omega}h^{\omega} \wedge
\widetilde{h}^{\omega} \widetilde{h}^{\omega}  -
\widetilde{h}^{\omega} \otimes  h^{\omega}h^{\omega} \wedge
(\widetilde{h}^{\omega}  h^{\omega} + h^{\omega}
\widetilde{h}^{\omega} )  \in {\sl ker}\, s \com (Id \otimes
\omega_{23}) \cong {\sl H}, \\
 h^{\omega} \otimes h^{\omega}h^{\omega} \wedge
(\widetilde{h}^{\omega} h^{\omega}  + h^{\omega}
\widetilde{h}^{\omega}) \in {\sl ker}\, \omega_{12}  \com (Id
\otimes \omega_{23}) \cong {\sl S}^3 {\sl H}.
\end{array}
$$
For $x,y,z \in {\sl E}$, we regard the tensor
$$
 \eta = \{ x^{\omega} \otimes( y^{\omega} \wedge z^{\omega} )
\wedge \omega_{\sl E}
 \} \otimes \left\{  h^{\omega} \otimes h^{\omega}h^{\omega}  \wedge
\widetilde{h}^{\omega} \widetilde{h}^{\omega}  -
\widetilde{h}^{\omega} \otimes  h^{\omega}h^{\omega} \wedge
(\widetilde{h}^{\omega}  h^{\omega} + h^{\omega}
\widetilde{h}^{\omega} )  \right\}.
$$
It is obvious that $ \eta  \in   {\sl E} \otimes \Lambda^2 {\sl E}
\wedge \{ \omega_{\sl E} \} \otimes {\sl ker} s \com (I \otimes
\omega_{23}) \subseteq {\sl EH} \otimes ({\sl EH})^4$.

 Now, if we consider the ${\sl Sp}(n){\sl Sp}(1)$-map $
Id\otimes {\sl a} : {\sl EH} \otimes  ({\sl EH})^4
 \to {\sl EH} \otimes  \Lambda^4 ({\sl EH}),
$ where $Id$ is the identity on ${\sl EH}$ and ${\sl a}$ is the
alternation map, we will obtain {\footnotesize
\begin{eqnarray*}
(Id \otimes {\sl a}) (\eta) & = & x^{\omega} h^{\omega} \otimes
y^{\omega} h^{\omega} \wedge z^{\omega} h^{\omega} \wedge
e^{\omega}_l \widetilde{h}^{\omega} \wedge
\widetilde{e}^{\omega}_l \widetilde{h}^{\omega} - x^{\omega}
h^{\omega} \otimes y^{\omega} \widetilde{h}^{\omega} \wedge
z^{\omega} \widetilde{h}^{\omega} \wedge e^{\omega}_l
h^{\omega} \wedge \widetilde{e}^{\omega}_l h^{\omega} \\
& & - x^{\omega} \widetilde{h}^{\omega} \otimes y^{\omega}
h^{\omega} \wedge z^{\omega} h^{\omega} \wedge e^{\omega}_l
\widetilde{h}^{\omega} \wedge \widetilde{e}^{\omega}_l h^{\omega}
 - x^{\omega} \widetilde{h}^{\omega} \otimes y^{\omega} h^{\omega}
\wedge z^{\omega}h^{\omega} \wedge e^{\omega}_l
h^{\omega} \wedge \widetilde{e}^{\omega}_l \widetilde{h}^{\omega} \\
&& + x^{\omega} \widetilde{h}^{\omega} \otimes y^{\omega}
\widetilde{h}^{\omega} \wedge z^{\omega} h^{\omega} \wedge
e^{\omega}_l h^{\omega} \wedge \widetilde{e}^{\omega}_l h^{\omega}
+ x^{\omega} \widetilde{h}^{\omega} \otimes y^{\omega} h^{\omega}
\wedge z^{\omega} \widetilde{h}^{\omega} \wedge e^{\omega}_l
h^{\omega} \wedge \widetilde{e}^{\omega}_l h^{\omega}.
\end{eqnarray*} }\noindent
Using identities \ref{equiache}, we can obtain  the restriction
$(Id \otimes {\sl a}) (\eta)_{|{\cal V}}$ of $(Id \otimes {\sl a})
(\eta)$ to ${\cal V}^5 \cong \left( [{\sl EH}] \right)^5 $.
Furthermore,  the real part ${\sl Re}(Id \otimes {\sl a})
(\eta)_{|{\cal V}}$ of this complex valued map  belongs to ${\cal
V}^{\ast} \otimes \Lambda^4 {\cal V}^{\ast}$ and, by an
straightforward computation, one can check  $ {\cal L}\left( {\sl
Re}(Id \otimes {\sl a}) (\eta)_{|{\cal V}} \right)= 4  {\sl Re}(Id
\otimes {\sl a}) (\eta)_{|{\cal V}}. $

If we choose $x,y,z \in {\sl E}$ such that $x \otimes y \wedge z$
belongs successively   to $\Lambda_0^3 {\sl E}$, ${\sl K}$ and
${\sl E}$, then ${\sl Re}(Id \otimes {\sl a}) (\eta)_{|{\cal V}}$
is an element of the irreducible real ${\sl Sp}(n){\sl
Sp}(1)$-representation determined by $\Lambda_0^3 {\sl EH}$, ${\sl
KH}$ and ${\sl EH}$, respectively. Therefore, by Schur's Lemma,
 any element $a$ of these ${\sl Sp}(n){\sl Sp}(1)$-modules
satisfies ${\cal L}(a) = 4 a$.
\vspace{2mm}

On the other hand, if we consider,  for $x,y,z \in {\sl E}$, the
tensor
$$
 \zeta = \{ x^{\omega} \otimes( y^{\omega} \wedge z^{\omega} )
\wedge \omega_{\sl E}  \} \otimes \left\{ h^{\omega} \otimes
h^{\omega}h^{\omega}  \wedge  (\widetilde{h}^{\omega} h^{\omega} +
h^{\omega} \widetilde{h}^{\omega}) \right\},
$$
we have
$
\zeta  \in  \left( {\sl E} \otimes \Lambda^2 {\sl E} \wedge \{
\omega_{\sl E} \} \right)  \otimes {\sl ker}  \omega_{12} \com (I
\otimes \omega_{23}) \subseteq {\sl EH} \otimes
({\sl EH})^4$.
 Therefore, {\footnotesize
$$
\begin{array}{l}
(Id \otimes {\sl a}) (\zeta)  =  2 x^{\omega} h^{\omega} \otimes
y^{\omega} h^{\omega} \wedge z^{\omega} h^{\omega} \wedge
e^{\omega}_l \widetilde{h}^{\omega} \wedge
\widetilde{e}^{\omega}_l h^{\omega} + 2 x^{\omega} h^{\omega}
\otimes y^{\omega} h^{\omega} \wedge z^{\omega} h^{\omega} \wedge
e^{\omega}_l
h^{\omega} \wedge \widetilde{e}^{\omega}_l \widetilde{h}^{\omega} \\
\qquad \qquad \quad - 2 x^{\omega} h^{\omega} \otimes y^{\omega}
\widetilde{h}^{\omega} \wedge z^{\omega} h^{\omega} \wedge e^{\omega}_l
h^{\omega} \wedge \widetilde{e}^{\omega}_l h^{\omega}
 - 2 x^{\omega} h^{\omega} \otimes y^{\omega} h^{\omega} \wedge
z^{\omega} \widetilde{h}^{\omega} \wedge e^{\omega}_l
h^{\omega} \wedge \widetilde{e}^{\omega}_l h^{\omega}.
\end{array}$$}\noindent
Taking equations \ref{equiache} into account, we obtain the
restriction $(Id \otimes {\sl a}) (\zeta)_{|{\cal V}}$ of $(Id
\otimes {\sl a}) (\zeta)$ to ${\cal V}^5 \cong \left( [EH]
\right)^5 $. The real part ${\sl Re}(Id \otimes {\sl a})
(\zeta)_{|{\cal V}}$ of this complex valued map  is in ${\cal
V}^{\ast} \otimes \Lambda^4 {\cal V}^{\ast}$. As before, one can
check $ {\cal L}\left( {\sl Re}(Id \otimes {\sl a})
(\zeta)_{|{\cal V}} \right)= -2   {\sl Re}(Id \otimes {\sl a})
(\zeta)_{|{\cal V}}. $

If we choose $x,y,z \in {\sl E}$ such that $x \otimes y \wedge z$
belongs successively   to $\Lambda_0^3 {\sl E}$, ${\sl K}$ and
${\sl E}$, then ${\sl Re}(Id \otimes {\sl a}) (\zeta)_{|{\cal V}}$
is an element of the irreducible real ${\sl Sp}(n){\sl
Sp}(1)$-representation determined by $\Lambda_0^3 {\sl E}{\sl S}^3
{\sl H}$, ${\sl K}{\sl S}^3 {\sl H}$ and ${\sl E}{\sl S}^3 {\sl
H}$, respectively. Therefore, by Schur's Lemma,
 any element $a$ of these ${\sl Sp}(n){\sl Sp}(1)$-modules
satisfies ${\cal L}(a) = -2 a$.

Since ${\cal L}$ is bijective  on each irreducible ${\sl
Sp}(n){\sl Sp}(1)$-module of ${\sl EH} \otimes \Lambda_0^2 {\sl E}
{\sl S}^2 {\sl H}$, then ${\cal L}$ is an ${\sl Sp}(n){\sl
Sp}(1)$-isomorphism. \cqd \vspace{2mm}

\seccion{ ${\sl Sp}(n){\sl Sp}(1)$-spaces of skew-symmetric
three-forms  }{\indent} \setcounter{proposition}{0} In this
section we study the irreducible modules which take part in the
${\sl Sp}(n){\sl Sp}(1)$-decomposition of $\Lambda^3 {\cal
V}^{\ast} \cong \Lambda^3 (EH) $,  given by $ \Lambda^3 {\cal
V}^{\ast} = ({\sl K} + {\sl E}){\sl H} + (\Lambda_{0}^3 {\sl E}  +
{\sl E}) {\sl S}^3 {\sl H}$(\cite{Sw2}). For such a purpose, we
consider the ${\sl Sp}(n){\sl Sp}(1)$-map ${\sl L}$ on $\Lambda^3
{\cal V}^{\ast}$  defined by \ref{lmant}. Relative to this ${\sl
Sp}(n){\sl Sp}(1)$-map, we have the following results.
\begin{proposition} \label{bb}
{\it
\begin{enumerate}
\item[{\rm i)}] ${\sl L}$ is an ${\sl Sp}(n){\sl
Sp}(1)$-isomorphism such that ${\sl L}^2 = 9 Id$.
 \item[{\rm ii)}]
$ ({\sl K} + {\sl E}){\sl H}$ consists of those three-forms $ b
\in \Lambda^{3} {\cal V}^{\ast}$ such that ${\sl L}(b) = 3 b$.
\item[{\rm iii)}] $(\Lambda_0^{2} {\sl E} + {\sl E}) {\sl S}^3
{\sl H}$  consists of those three-forms $b \in \Lambda^3 {\cal
V}^{\ast}$  such that ${\sl L}(b) = -3 b$.
\end{enumerate}
}
\end{proposition}
\Dt Part i) follows by an straightforward computation. Now, we
consider the complex tensor
$
\eta = {\sl a} \left(  (x^{\omega} \otimes y^{\omega} \wedge
z^{\omega})  \otimes \left( 2
\widetilde{h}^{\omega} \otimes h^{\omega} h^{\omega} - h^{\omega}
\otimes (\widetilde{h}^{\omega} h^{\omega} + h^{\omega}
\widetilde{h}^{\omega}  \right)  \right),
$
where ${\sl a} \, : \, ({\sl EH})^3 \to \Lambda^3 ({\sl EH})$ is
the alternation map.

 Since
$
2 \widetilde{h}^{\omega} \otimes h^{\omega} h^{\omega} -
h^{\omega} \otimes (\widetilde{h}^{\omega} h^{\omega} + h^{\omega}
\widetilde{h}^{\omega}  ) \in {\sl ker}\,{\sl s}\cong {\sl H},
$
where  ${\sl s}
 : {\sl H} \otimes {\sl S}^2 {\sl H} \to {\sl S}^3 {\sl H}$  is the symmetrization map,
 then $\eta \in ( {\sl K} + {\sl E}) {\sl  H} \subseteq \Lambda^3
({\sl EH})$.

The real part ${\sl Re}(\eta)$ of the restriction of $\eta$ to
$[{\sl EH}]$ is an element  of the real ${\sl Sp}(n){\sl
Sp}(1)$-module determined by $({\sl K} + {\sl E}) {\sl H}$
($(+1)$-eigenspace of the real structure map), i.e.,
{\footnotesize
\begin{eqnarray*}
{\sl Re}(\eta) & = &  -2 J x \wedge y \wedge z - 2 K x \wedge I y
\wedge  z  -  2 K x \wedge y \wedge I z + 2
J x  \wedge I y \wedge I z  \\
&  &  + x \wedge J y \wedge z + I x  \wedge K y \wedge z -  I x
\wedge J y \wedge I z +
 x  \wedge K y \wedge  I z  \\
&  &  + x \wedge  y \wedge J z - I x  \wedge I y \wedge J z + I x
\wedge  y \wedge K z +
 x  \wedge I y \wedge  K z
\end{eqnarray*}
} is contained in $[{\sl K} {\sl H}] + [{\sl  EH}] \subseteq
\Lambda^3 [{\sl EH}]$. Now, by straightforward computation we get
${\sl L}({\sl Re}(\eta)) = 3 {\sl Re}(\eta)$. Hence we have part
ii). \vspace{2mm}

Finally, we consider $ \zeta = {\sl a} \left(  (x^{\omega} \otimes
y^{\omega} \wedge z^{\omega} )  \otimes ( h \otimes h h ) \right).
$ Taking $h\otimes hh \in {\sl ker}\, \omega_{12} \cong {\sl S}^3
{\sl H}$ into account, where $\omega_{12} \, : {\sl H} \otimes
{\sl S}^2 {\sl H} \to {\sl H}$,  we have $\zeta \in (\Lambda_{0}^3
{\sl E}  + {\sl E}) {\sl S}^3 {\sl H}$. The real part ${\sl
Re}(\zeta)$ of the restriction of $\zeta$ to $[EH]$ is an element
of the real ${\sl Sp}(n){\sl Sp}(1)$-module determined by $
(\Lambda_{0}^3 {\sl E} + {\sl E}) {\sl S}^3 {\sl H} $, i.e.,
{\footnotesize
\begin{eqnarray*}
{\sl Re}(\zeta) & = &   x \wedge y \wedge z - I x \wedge I y
\wedge  z  -  I x \wedge y \wedge I z -
 x  \wedge I y \wedge I z
\end{eqnarray*} }
belongs to $[\Lambda_{0}^3 {\sl E} {\sl S}^3 {\sl H}] + [{\sl E}
{\sl S}^3 {\sl H}] \subseteq \Lambda^3 [{\sl EH}]$. By
 straightforward computation we get  ${\sl L}({\sl Re}(\zeta)) =
- 3 {\sl Re}(\zeta)$. Thus part iii) follows. \cqd
\begin{observation} {\it
 Note that  the condition  ${\sl L}(b)=3b$ is equivalent to
$(A_{(1)} A_{(2)} + A_{(2)} A_{(3)} + A_{(3)} A_{(1)}) b =b$, for
$A=I,J,K$ and, on the other hand,  the condition  ${\sl L}(b)=-
3b$ is equivalent to $\sum_{A=I,J,K} A_{(2)} A_{(3)} b = -b$.
Therefore, we also  have these alternative ways to describe the
spaces $ ({\sl K} + {\sl E}){\sl H}$ and $(\Lambda_0^{2} {\sl E} +
{\sl E}) {\sl S}^3 {\sl H}$.  }
\end{observation}


The  ${\sl Sp}(n){\sl Sp}(1)$-subspace ${\sl EH}$ of $\Lambda^3
{\cal V}^{\ast}$ consists of  the three-forms $x \lrcorner
\Omega$, where $x \in {\cal V}$. For $b \in \Lambda^3 {\cal
V}^{\ast}$, its projection on ${\sl EH}$ is given by $\pi_{\sl
EH}(b) = \xi_b \lrcorner \Omega$, where $\xi_b$ is defined by
\begin{equation} \label{xi}
\qquad \xi_b(x)  = \mbox{$\frac{1}{12 (2n+1)}$}  \sum_{A=I,J,K}
b(e_r , A e_r , A x )= - \mbox{$\frac{1}{6(2n+1)}$} \sum_{A=I,J,K}
\langle Ax \lrcorner b, \omega_A \rangle,
\end{equation}
for all $x \in {\cal V}$, where $\{ e_1, \dots , e_{4n} \}$ is an
orthonormal basis for ${\cal V}$. Note that the map $\Lambda^{3}
{\cal V}^{\ast} \to {\cal V}^{\ast}$, given by $b \to \xi_b$,
 is an ${\sl Sp}(n){\sl Sp}(1)$-map.
\vspace{2mm}

The subspace of $\Lambda^3 {\cal V}^{\ast}$ consisting of the
three-forms  $ b = - 2 \sum_{A=I,J,K} A \xi_A \wedge \omega_{A}$,
where $\xi_I, \xi_J , \xi_K \in {\cal V}^{\ast}$, is an ${\sl
Sp}(n){\sl Sp}(1)$-module of dimension $12n$. By reasons of
dimension,   such an ${\sl Sp}(n){\sl Sp}(1)$-module must coincide
with ${\sl E}( {\sl  H} + {\sl S}^3 H)$. Moreover, if $\xi_b
\lrcorner \Omega$ is the ${\sl EH}$-projection of  $ b = - 2
\sum_{A=I,J,K}A \xi_{A} \wedge \omega_{A}$, then $ \xi_b =
\frac{1}{3}( \xi_I + \xi_J + \xi_K)$.

Therefore, the subspace ${\sl E}{\sl S}^3 {\sl H}$ of $\Lambda^3
{\cal V}^{\ast}$ consists of the three-forms $ b = \linebreak - 2
\sum_{A=I,J,K} A \xi_{A} \wedge \omega_{A}$ such that $\xi_I +
\xi_J + \xi_K=0$. That is, the kernel of the map ${\sl E}( {\sl H}
+{\sl S}^3 {\sl H}) \to {\sl EH}$ given by $b \to \xi_b$. Hence
the subspace ${\sl EH}$ can be also described as consisting of
those three-forms $ b = - 2 \sum_{A=I,J,K} A \xi_{A} \wedge
\omega_{A}$ such that $\xi_I = \xi_J = \xi_K$. Note that, in such
a situation, $\xi_b = \xi_{I}=\xi_J = \xi_K$ and $b = \xi_b
\lrcorner \Omega$.

For $b \in \Lambda^3 {\cal V}^{\ast}$, $  \pi_{{\sl E}({\sl H} +
{\sl S}^3 {\sl H})}(b) = - 2 \sum_{A=I,J,K} A \xi_{b;A} \wedge
\omega_{A}$ is its projection on ${\sl E} ({\sl H} +{\sl S}^3{\sl
H}) $, where
\begin{equation} \label{xialpha}
\xi_{b;A} (x) = - \mbox{$\frac{3}{2(n-1)}$} \xi_b (x)  -
\mbox{$\frac{1}{4(n-1)}$}
 \langle Ax \lrcorner b, \omega_A \rangle,
\end{equation}
for all $x \in {\cal V}$. Note that each one-form $\xi_{b;A}$ is
not global, because its definition depends on the chosen basis
$\{I,J,K\}$. However, the three-form $-  2 \sum_{A=I,J,K} A
\xi_{b;A} \wedge \omega_{A}$ is global and  $b \to -  2
\sum_{A=I,J,K} A \xi_{b;A} \wedge \omega_{A} $ constitutes an
${\sl Sp}(n){\sl Sp}(1)$-map.

After these considerations and taking Proposition \ref{bb} into
account, in Table 1 we display the corresponding conditions which
describe the ${\sl Sp}(n){\sl Sp}(1)$-subspaces of $\Lambda^{3}
{\cal V}^{\ast} \cong \Lambda^3 ({\sl EH})$. In case of ${\cal V}$
 eight-dimensional, we would have $\Lambda_{0}^3 {\sl E} {\sl
S}^3 {\sl H} = \{ 0 \}$ and a smaller corresponding table would be
obtained.

{\scriptsize
\begin{center}
{\normalsize  Table 1: ${\sl Sp}(n){\sl Sp}(1)$-subspaces of
$\Lambda^3 {\cal V}^{\ast}\cong \Lambda^3 ({\sl EH})$}
\vspace{1mm}

\begin{tabular}{|p{3.3cm}|p{6.7cm}|}
\hline
 &\\[-3mm]
$\{ 0 \}$ & $b =0$ \\[1mm]
\hline
 &\\[-3mm]
{\scriptsize ${\sl K} {\sl H}$}& $ \displaystyle {\sl L}(b) = 3 b$
or $\displaystyle \sum_{1\leq i < j \leq 3}A_{(i)} A_{(j)} b =b$,
for $A=I,J,K$,  and $\xi_b = 0 $
\\[1mm]
\hline
 &\\[-3mm]
{\scriptsize ${\sl E} {\sl H}$}& $b = \xi_b \lrcorner \Omega
$    \\[1mm]
\hline
 &\\[-3mm]
{\scriptsize $\Lambda_0^3 {\sl E} {\sl S}^3 {\sl H}$}&
 $ {\sl L}(b) =  - 3 b$  or $\displaystyle \sum_{A=I,J,K} A_{(2)} A_{(3)} b =-b$,
\newline   and $  \xi_{b;I}  = \xi_{b;J} = \xi_{b,K}= 0$
 \\[1mm]
\hline
 &\\[-3mm]
{\scriptsize $ {\sl E} {\sl S}^3  {\sl H} $}&
 $ b  = - 2 \displaystyle \sum_{A=I,J,K} A \xi_{b;A}  \wedge
\omega_{A} $ and $\xi_b  = 0$    \\[1mm]
\hline
 &\\[-3mm]
{\scriptsize $ ({\sl K} +  {\sl E})H $}& ${\sl L}(b)=  3 b$ or
$\displaystyle \sum_{1\leq i < j \leq 3}A_{(i)} A_{(j)} b =b$, for
$A=I,J,K$
 \\[1mm]
\hline
 &\\[-3mm]
{\scriptsize $ {\sl KH} +  \Lambda_{0}^3 {\sl E} {\sl S}^3 {\sl H}
$}& $\xi_{b;I} = \xi_{b;J} = \xi_{b;K}
=0$  \\[1mm]
\hline
 &\\[-3mm]
{\scriptsize $ {\sl K H} + {\sl  E} {\sl S}^3  {\sl H} $}& $ {\sl
L}(b) = 3 b + 12 \displaystyle \sum_{A=I,J,K}  \xi_{b;A} \wedge
\omega_{A}$ \\[1mm]
\hline
 &\\[-3mm]
{\scriptsize $    {\sl  EH} + \Lambda^3_{0}  {\sl E}  {\sl S}^3
{\sl H}$}&  ${\sl L}(b)  = - 3b + 6 \xi_b \lrcorner \Omega $ and
$\xi_{b;I} = \xi_{b;J} = \xi_{b;K}$
 \\[1mm]
\hline
 &\\[-3mm]
{\scriptsize $ {\sl E} ( {\sl  H} + {\sl S}^3  {\sl H}) $}& $ b =
- 2 \displaystyle \sum_{A=I,J,K} A\xi_{b;A}
\wedge \omega_{A}$  \\[1mm]
\hline
\end{tabular}
\end{center}
}

\newpage
{\scriptsize
\begin{center}
\begin{tabular}{|p{3.3cm}|p{6.7cm}|}
\hline
 &\\[-3mm]
{\scriptsize $(\Lambda_0^3 {\sl E} + {\sl E}) {\sl S}^3 {\sl H}$}&
$
{\sl L}(b) = - 3b $ or  $\displaystyle \sum_{A=I,J,K} A_{(2)} A_{(3)} b =-b$  \\[1mm]
\hline
 &\\[-3mm]
{\scriptsize $ ({\sl K} + {\sl E})  H +  \Lambda_0^3 {\sl E} {\sl
S}^3
 {\sl H} $}&
 $ \xi_{b;I} = \xi_{b;J} = \xi_{b;K}$ \\[1mm]
\hline
 &\\[-3mm]
{\scriptsize $ ({\sl K} + {\sl E}){\sl H} + {\sl E}  {\sl S}^3
{\sl H} $}& ${\sl L} (b) = 3 b + 6 \xi_b \lrcorner \Omega +12
\displaystyle
\sum_{A=I,J,K} A \xi_{b;A} \wedge \omega_{A}$  \\[1mm]
\hline
 &\\[-3mm]
{\scriptsize ${\sl K} {\sl H} + (\Lambda^3_0  {\sl E} + {\sl E}) {\sl S}^3 H$}&  $ \xi_b = 0$ \\[1mm]
\hline
&\\[-3mm]
{\scriptsize $ {\sl E} {\sl H} + (\Lambda_0^3 {\sl E} + {\sl E})
{\sl S}^3 {\sl H}$}&
 $ {\sl L}(b) = - 3b + 6  \xi_b \lrcorner \Omega $   \\[1mm]
\hline
 &\\[-3mm]
{\scriptsize $({\sl K} + {\sl E}) {\sl H} + ( \Lambda_0^3 {\sl E}
+ {\sl E}) {\sl S}^3 {\sl H}$}&
no relation   \\[1mm]
\hline
\end{tabular}
\end{center}
}
\vspace{1mm}


\seccion{Classes of almost quaternion-Hermitian manifolds
}{\indent} \setcounter{proposition}{0} In this  section, our aim
is to deduce  conditions which characterize the different classes
of almost quaternion-Hermitian manifolds. In other words, we are
going  to show the explicit conditions  which describe the ${\sl
Sp}(n){\sl Sp}(1)$-subspaces of the space ${\cal V}^{\ast} \otimes
\Lambda_0^2 E {\sl S}^2 H$ of covariant derivatives of the
fundamental four-form $\Omega$.

For such a purpose,  we consider the ${\sl Sp}(n){\sl Sp}(1)$-map
$ {\sl d}^*  \, : \, {\cal V}^{\ast}   \otimes \Lambda_0^2 {\sl E}
{\sl S}^2 {\sl H}  \to \Lambda^3 {\cal V}^{\ast}, $ defined by
${\sl d}^* a =  - {\sl C}_{12}(a) $, where ${\sl C}$ is  the
metric contraction. In \cite{Sw4} it is proved
\begin{equation} \label{swdelta}
{\sl ker} \, {\sl d}^* = \Lambda_0^3 {\sl E H} + {\sl K} {\sl S}^3
{\sl H}.
\end{equation}

If $b$ is a three-form, one can  check that the tensors ${\sl
i}_{A} (\cdot \lrcorner b)$ satisfy i)$'$ and ii)$''$ of
Proposition \ref{Perpsp}. Then $ \sum_{A=I,J,K} {\sl i}_{A} (
\cdot \lrcorner b ) \wedge \omega_{A} \in {\cal V}^{\ast} \otimes
\Lambda_0^2 {\sl E} {\sl S}^2 {\sl H}$. Therefore, we can consider
the  ${\sl Sp}(n){\sl Sp}(1)$-map $ \widehat{\sl d}^* \, : \,
\Lambda^3 {\cal V}^{\ast} \to ({\sl K} + {\sl E}) {\sl H} +
(\Lambda_0^3 {\sl E} + {\sl E}) {\sl S}^3 {\sl H}, $ defined by
\begin{eqnarray*}
\widehat{\sl d}^*(b) & = &  \mbox{$\frac{1}{18}$} \sum_{A=I,J,K}
{\sl i}_{A} (\cdot \lrcorner {\sl L}(b)) \wedge \omega_{A}  -
\mbox{$\frac{2k_1}{3k_2}$} \sum_{A,B=I,J,K} {\sl i}_{A} (\cdot
\lrcorner (B\xi_{b;B} \wedge \omega_{B} ) ) \wedge \omega_{A} \\
 && - \mbox{$\frac{4 k_1^2 + k_2^2}{12 k_1k_2}$} \{ \cdot \wedge
  (\xi_b  \lrcorner \Omega ) - \xi_b \wedge (\cdot
\lrcorner \Omega) \},
\end{eqnarray*}
 where $k_1=n-1$, $k_2=2n+1$ and   ${\sl L}$, $\xi_b$, $\xi_{b;A}$ are defined by
 \ref{lmant}, \ref{xi}, \ref{xialpha}, respectively. By straightforward
computation one can check ${\sl d}^* \com \widehat{\sl d}^* (b) =
b$. Therefore, ${\sl d}^*$ is surjective and $\widehat{\sl d}^*$
is injective.

Next we will show a relation involving the maps ${\sl L}$, ${\sl
d}^*$ and ${\sl L}$.
\begin{lemma}
{\it Let  ${\cal L}$ and ${\sl L}$ be the maps defined by
\ref{eleauto} and \ref{lmant}, respectively. Then
\begin{equation} \label{ldeltab}
{\sl d}^* \com {\cal L} = {\sl d}^* + {\sl L} \com {\sl d}^*.
\end{equation}
}
\end{lemma}
\Dt In fact, let  $a = a_1 + a_2 \in {\cal V}^{\ast} \otimes
\Lambda_0^2 {\sl E} {\sl S}^2 {\sl H}$ , where $a_1 \in
(\Lambda_0^3 {\sl E} + {\sl K} + {\sl E}){\sl H}$ and $a_2 \in
(\Lambda_0^3 {\sl E} + {\sl K} + {\sl E}) {\sl S}^3 {\sl H}$.
Taking Proposition \ref{laele} into account, we have $ {\sl d}^*
{\cal L}(a) = 4 {\sl d}^* a_1 - 2 {\sl d}^* a_2. $ Since ${\sl
d}^* a_1$ and ${\sl d}^* a_2$ are the projections of ${\sl d}^* a$
on $({\sl K} + {\sl E}){\sl H}$ and $(\Lambda_0^3 {\sl E} + {\sl
E}){\sl S}^3 {\sl H}$, respectively, we have ${\sl d}^* a_1 =
\frac{1}{6} (3 {\sl d}^* a +{\sl L}({\sl d}^* a))$ and ${\sl d}^*
a_2 = \frac{1}{6}(3 {\sl d}^* a - {\sl L}({\sl d}^* a))$ and the
required identity holds. \cqd \vspace{2mm}

By reiterated use of the maps  ${\cal L}$, $\widehat{\sl d}^*$ and
Table 1,  taking Schur's Lemma into account,  we can get the
conditions which define each ${\sl Sp}(n){\sl Sp}(1)$-subspace of
${\cal V}^{\ast} \otimes \Lambda_0^2 E {\sl S}^2 H$.
 Identity \ref{ldeltab} is used to eliminate some redundant
conditions. Thus, the different classes are characterized by
conditions displayed in the second column of Table 2, where we
denote $k_1=n-1$, $k_2=2n+1$, $\xi = \xi_{{\sl d}^* \Omega}$ and
$\xi_{A} = \xi_{{\sl d}^* \Omega;A}$.

We will show some examples to illustrate the  method we follow.
\vspace{2mm}

\noindent {\sc EXAMPLE 1}: Let $ \nabla \Omega \in {\sl KH}$, then
$\nabla \Omega = \widehat{\sl d}^* ({\sl d}^* \Omega)$ (we are
writing ${\sl d}^* \Omega = {\sl d}^* (\nabla \Omega)$, in the
context of manifolds ${\sl d}^*  \Omega$ means the coderivative of
$\Omega$). From Table 1, we have ${\sl L}({\sl d}^* \Omega ) = 3
{\sl d}^* \Omega$ and $\xi = 0$, then
\begin{equation} \label{KH}
\nabla_x \Omega = \mbox{$\frac{1}{6}$}  \sum_{A=I,J,K}  {\sl
i}_{A} ( x \lrcorner {\sl d}^* \Omega )  \wedge \omega_{A}.
\end{equation}
for all $x \in {\cal V}$.
Thus, we could say:
\vspace{1mm}

 $\nabla \Omega  \in {\sl KH}$ if and only if
 $\nabla \Omega$ satisfies \ref{KH}, ${\sl L}({\sl d}^* \Omega) = 3
{\sl d}^* \Omega$ and $\xi=0$. \vspace{1mm}

But some of these conditions are redundant. In fact, from the
equation \ref{KH},
taking  \ref{xi} and \ref{xialpha} into account,  we have
\begin{equation} \label{KH1}
{\sl d}^* \Omega = \mbox{$\frac{1}{3}$} {\sl L}({\sl d}^* \Omega)
- \xi \lrcorner \Omega + \mbox{$\frac{4(n-1)}{3}$} \sum_{A=I,J,K}
A \xi_{A} \wedge \omega_{A}.
\end{equation}
Since ${\sl L}^2({\sl d}^* \Omega) = 9 {\sl d}^* \Omega$, ${\sl
L}(\xi \lrcorner \Omega) = 3 \xi \lrcorner \Omega$ and ${\sl L}
\left(  \sum_{A=I,J,K} A(\xi_{A} - \xi) \wedge \omega_{A} \right)
= \linebreak - 3  \sum_{A=I,J,K} A (\xi_{A} - \xi)  \wedge
\omega_{A}$, we have
\begin{equation} \label{KH2}
{\sl L}({\sl d}^*  \Omega) = 3 {\sl d}^* \Omega - (4n-1) \xi
\lrcorner \Omega - 4 (n-1) \sum_{A=I,J,K} A \xi_{A} \wedge
\omega_{A}.
\end{equation}
Now, equations \ref{KH1} and \ref{KH2} imply  $\xi \lrcorner
\Omega = 0$. Therefore, $\xi = 0$.

Taking $\xi = \frac{1}{3} (\xi_I + \xi_J + \xi_K)$ and equation
\ref{KH1} into account,  it follows next description of the
$KH$-subspace:

$\nabla \Omega  \in {\sl KH}$ if and only if $\nabla \Omega$
satisfies \ref{KH} and $\xi_I = \xi_J = \xi_K$. \vspace{2mm}

\noindent {\sc EXAMPLE 2}: Let $ \nabla \Omega \in \Lambda_0^3
{\sl E} ( {\sl H} + {\sl S}^3 {\sl H})$, then $\nabla \Omega =
(\nabla \Omega)_1 + (\nabla \Omega)_2$, where $(\nabla \Omega)_1
\in \Lambda_0^3 {\sl EH}$ and $(\nabla \Omega)_2 \in \Lambda_0^3
{\sl E} {\sl S}^3 {\sl H}$. By \ref{swdelta}, ${\sl d}^* \Omega =
{\sl d}^* (\nabla \Omega)_2$. Since ${\sl d}^* \Omega \in
\Lambda_0^3 {\sl E}{\sl S}^3 {\sl H}$, then ${\sl L}({\sl d}^*
\Omega) = - 3 {\sl d}^* \Omega$ and $\xi_I = \xi_J = \xi_K$ (see
Table 1). Therefore, we have \vspace{-7mm}

$$
(\nabla_{\cdot} \Omega)_2 = \widehat{\sl d}^*  ({\sl d} \Omega) = -
\mbox{$\frac{1}{6}$} \sum_{A=I,J,K}  {\sl i}_{A}  ( \cdot
\lrcorner {\sl d}^* \Omega )  \wedge \omega_{A}.
$$\vspace{-5mm}

From Proposition \ref{laele} we obtain
\begin{equation} \label{LEHLES3H}
{\cal L}(\nabla_{\cdot} \Omega) = 4 (\nabla_{\cdot} \Omega )_1 - 2
(\nabla_{\cdot} \Omega)_2 = 4 \nabla_{\cdot} \Omega +
\sum_{A=I,J,K} {\sl i}_{A}  ( \cdot \lrcorner {\sl d}^* \Omega )
\wedge \omega_{A}.
\end{equation} \vspace{-6mm}

\noindent Thus, taking \ref{ldeltab} into account,  we could say:

$\nabla \Omega  \in \Lambda_0^3 {\sl E} ( {\sl H} + {\sl S}^3 {\sl
H}) $ if and only if
 $\nabla \Omega$ satisfies \ref{LEHLES3H},
 ${\sl L}({\sl d}^* \Omega) = - 3 {\sl d}^* \Omega$ and $\xi_I = \xi_J =
\xi_K$. \vspace{1mm}

But again, some of these conditions are redundant. In fact, from the
equality \ref{LEHLES3H},
taking   \ref{xi} and \ref{xialpha} into account,  we have
\begin{equation} \label{LEH1}
\qquad {\sl d}^* {\cal L}(\nabla \Omega)  =  {\sl d}^* \Omega +
{\sl L}({\sl d}^* \Omega) = 4 {\sl d}^* \Omega + 2  {\sl L}({\sl
d}^* \Omega) - 6  \xi \lrcorner \Omega + 8(n-1) \sum_{A=I,J,K} A
\xi_{A}  \wedge \omega_{A}.
\end{equation} \vspace{-5mm}

\noindent Applying the map ${\sl L}$ to both sides of this
equation, we get
\begin{equation} \label{LEH2}
0 = 3 {\sl d}^*  \Omega + {\sl L}({\sl d}^* \Omega)  - 2 (4n-1)
\xi \lrcorner \Omega - 8 (n-1) \sum_{A=I,J,K} A \xi_{A}  \wedge
\omega_{A}.
\end{equation}\vspace{-5mm}

\noindent From equations \ref{LEH1} and \ref{LEH2}, it follows   $
0 = 2 \sum_{A=I,J,K}A (\xi_{A} - \xi ) \wedge \omega_{A} $.
Therefore $\xi_I = \xi_J = \xi_K = \xi$.

Now, taking equation \ref{LEH1} into account, we obtain the
following description of the subspace $\Lambda_0^3 {\sl E} ({\sl
H} + {\sl S}^3 {\sl H})$:

$\nabla \Omega  \in \Lambda_0^3 {\sl E} ({\sl H} + {\sl S}^3 {\sl
H}) $ if and only if $\nabla \Omega$ satisfies \ref{LEHLES3H} and
$\xi = 0$. \vspace{1mm}

As we have already pointed out, by similar considerations to those
contained in the exposed examples,  we will get the conditions
given in Table 2 which describe how must be $\nabla \Omega$ to
belong to the different ${\sl Sp}(n){\sl Sp}(1)$-subspaces of
${\cal V}^{\ast} \otimes \Lambda_0^2 {\sl E} {\sl S}^2 {\sl H}$.

In the mentioned Table 2, the conditions are displayed for ${\cal
V} = {\sl T}_p M$, the tangent space at a point $p$ of an almost
quaternion-Hermitian manifold $M$ of dimension higher than eight.
 In case $M$ were eight-dimensional, we would have $\Lambda_0^3 {\sl
 E}
{\sl H} = \Lambda_0^3 {\sl E} {\sl S}^3 {\sl H} = \{ 0 \}$ and an
smaller corresponding table would be resultant. \vspace{1mm}

  A. Swann in \cite{Sw1} showed that, for dimension higher than eight,
    $\nabla \Omega$ is totally
determined by the exterior derivative ${\sl d}\Omega$. Therefore,
our purpose now is to describe the different classes of almost
quaternion-Hermitian manifolds by conditions on ${\sl d}\Omega$.
For this aim, we consider we consider the alternation map ${\sl a}
\, : \, {\sl T}^{\ast} M \otimes \Lambda_0^2 {\sl E} {\sl S}^2
{\sl H} \to \Lambda^5 {\sl T}^{\ast} M $, defined by ${\sl
a}(a)(X,Y,Z,U,V) = \sumcic_{XYZUV} a(X,Y,Z,U,V)$, where $\sumcic$
denotes cyclic sum. For dimension $4n > 8$,  it was shown in
\cite{Sw1} that ${\sl a}$ is non-zero on each irreducible summand.
Thus, by Schur's Lemma, applying ${\sl a}$ to the conditions
already obtained in terms of $\nabla \Omega$, we will get
conditions in terms of ${\sl d}\Omega$ which characterize, for
dimension $4n >8$, the different classes of almost
quaternion-Hermitian manifolds. Proceeding in this way, taking
Lemma below into account, we obtain the  conditions  displayed in
the third column of Table 2. \vspace{1mm}

If the almost quaternion-Hermitian manifold is eight-dimensional,
we have  $\Lambda_{0}^3 {\sl E}= \{ 0 \}$ and $\Lambda^5 {\sl
T}^{\ast} M \cong \Lambda^3 {\sl T}^{\ast} M = ({\sl K} + {\sl E})
H + {\sl E} {\sl S}^3 {\sl H}$. Therefore,
 a partial set of classes can be characterized  via ${\sl d} \Omega$.
  Table 3 contains such characterizations.
\begin{lemma} \label{dddcond}
{\it Given any one-form $\zeta$ and any three-form $b$ on an
almost quaternion-Hermitian manifold $M$, and the alternation map
${\sl a} \,  : \, {\sl T}^{\ast}  M \otimes \Lambda_0^2 {\sl E}
{\sl S}^2 {\sl H} \to \Lambda^5 {\sl T}^{\ast}  M $, then
\begin{eqnarray} \label{ddd1}
{\sl a} \left( {\cal L} (\nabla \Omega) \right) + 2 {\sl d} \Omega
& = & {\sl L}({\sl d} \Omega ), \\
\label{ddd2} {\sl a} \left(  \cdot \wedge (\zeta  \lrcorner
\Omega) - \zeta \wedge ( \cdot \lrcorner  \Omega) \right) & = & 4
\zeta \wedge \Omega, \\
 \label{ddd3} {\sl a} \left( \sum_{A=I,J,K}  {\sl i}_{A} ( \cdot
\lrcorner b) \wedge \omega_{A} \right) & = &2 \sum_{A=I,J,K} {\sl
i}_{A} b \wedge \omega_{A}.
\end{eqnarray} }
\end{lemma}
\Dt Equations \ref{ddd3} and \ref{ddd2}  follow by straightforward
computation. Taking Lemma \ref{lemma:La2a} into account, equation
\ref{ddd1} also  follows by direct computation. \cqd \vspace{1mm}

The coderivative ${\sl d}^* \Omega$ and one-forms $\xi$, $\xi_I$,
$\xi_J$, $\xi_K$ are involved in some of the conditions contained
in the third column of Table 2. Therefore, it is natural question
wether such objects can be computed by using  the exterior
derivative ${\sl d}\Omega$. We will answer positively to such a
question. In fact, one can check \vspace{-7mm}

$$
\mbox{$\frac{(-1)^{n+1}}{(2n+1)!}$} \Omega^n  =   e_1 \wedge \dots
\wedge e_n \wedge Ie_1 \wedge \dots \wedge I e_n  \wedge J e_1
\wedge \dots \wedge J e_n \wedge K e_1 \wedge \dots \wedge K e_n.
$$ \vspace{-5mm}

\noindent Thus, taking as volume form ${\sl Vol} = \frac{
(-1)^{n+1}}{(2n+1)!} \Omega^n$, we obtain \vspace{-3mm}
\begin{equation}
{\sl d}^*  \Omega   = - \ast {\sl d} \ast \Omega =  \mbox{$\frac{
(-1)^{n}(n-1)}{(2n+1)!}$} \ast (\Omega^{n-2} \wedge {\sl
d}\Omega). \label{coddif}
\end{equation} \vspace{-6mm}

For a one-form  $a$ on $M$, it can be checked that
\begin{equation} \label{estre}
\ast \left( \ast(a \wedge \Omega) \wedge \Omega\right) = 12
(n-1)(2n+1) a.
\end{equation}

If we consider the ${\sl Sp}(n){\sl Sp}(1)$-map ${\sl T}^{\ast}_p
M \otimes \Lambda^2_0 {\sl E} {\sl S}^2 {\sl H} \to {\sl
T}^{\ast}_p M$ given by $ \nabla \Omega \to  \ast( \ast ({\sl a}
(\nabla \Omega)) \wedge \Omega)$, where ${\sl a}$ is the
alternation map above defined,  then we get
$$
 \ast \left( \ast( {\sl a} (\nabla \Omega) \right) \wedge \Omega) =  \ast\left( \ast({\sl a}
(\pi_{\sl EH} (\nabla \Omega))) \wedge \Omega \right) = -
\mbox{$\frac{1}{n-1}$} \ast \left( \ast (\xi \wedge \Omega) \wedge
\Omega \right),
$$
Now, from Table 2 we have
\begin{equation} \label{ehcomp}
\qquad \pi_{\sl EH}(\nabla_{\cdot} \Omega ) = -
\mbox{$\frac{1}{4(n-1)}$} \{ \cdot \wedge (\xi \lrcorner \Omega )
- \xi \wedge (\cdot \lrcorner \Omega) \}.
\end{equation}
 Taking   ${\sl d}\Omega = {\sl a}(\nabla \Omega)$, \ref{estre}
and \ref{ehcomp} into account, we obtain the identity
\begin{equation} \label{xiast}
\xi = \mbox{$\frac{-1}{12(2n+1)}$} \ast(\ast {\sl d}\Omega \wedge
\Omega).
\end{equation}

Now, if $\zeta_I$, $\zeta_J$ and $\zeta_K$ are three one-forms, it
can be checked \vspace{-2mm}

\begin{equation} \label{estre1}
\ast \left(   \mbox{$\sum_{B=I,J,K}$} \ast \left(   B \zeta_{B}
\wedge \omega_{B} \right) \wedge \omega_{A} \right) = 4 n
A\zeta_{A}.
\end{equation}\vspace{-2mm}

If we consider the ${\sl Sp}(n){\sl Sp}(1)$-map $\Lambda^3 {\sl
T}^{\ast}_p M  \to {\sl E H} + {\sl E}{\sl S}^3 {\sl H}$ given by
$ {\sl d}^* \Omega \to \frac{1}{4n}  \sum_{A=I,J,K} \ast \left(
\ast {\sl d}^*   \Omega \wedge \omega_{A} \right) \wedge
\omega_{A}$, then, taking   identity \ref{estre1} into account, we
have \vspace{-7mm}

\begin{eqnarray*}
    \sum_{A=I,J,K} \ast \left( \ast ({\sl
d}^* \Omega) \wedge \omega_{A} \right) \wedge \omega_{A} &  = &
\displaystyle  \sum_{A=I,J,K} \ast \left( \ast(\pi_{{\sl E} ({\sl
H} + {\sl S}^3{\sl H})} {\sl d}^* \Omega) \wedge \omega_{A}
\right) \wedge \omega_{A} \\
& = &  - 8n \sum_{A=I,J,K} A\xi_{A} \wedge \omega_A.
\end{eqnarray*}\vspace{-6mm}

\noindent Hence $ A \xi_{A} = - \frac{1}{8n} \ast \left( \ast (
{\sl d}^* \Omega) \wedge \omega_{A} \right)$.  Therefore, by
\ref{coddif}, $\xi_A$ can be computed via ${\sl d} \Omega$.
\vspace{1mm}

 Now we
deduce alternative ways to characterize some classes of almost
quaternion-Hermitian manifolds. For such a purpose, we consider
the ${\sl Sp}(n){\sl Sp}(1)$-module of $\Lambda^5 T^{\ast} M$
given by ${\sl E} ({\sl H} + {\sl S}^3 {\sl H} + {\sl S}^5 {\sl
H})$ (see \cite{Sw2}), which consists of those five-forms
$\varphi$ such that $ \varphi = \sum_{A,B =I,J,K} a_{AB} \wedge
\omega_{A} \wedge \omega_{B}, $ where $a_{AB}$ are one-forms.
Therefore,
 the orthogonal complement  $\left( {\sl E} ( {\sl H} + {\sl
S}^3 {\sl H}+  {\sl S}^5 {\sl H}) \right)^{\perp}$ consists of
those five-forms $\varphi$ such that $ \langle a \wedge \omega_{A}
\wedge \omega_{B} , \varphi \rangle = 0, $ for $A,B=I,J,K$ and for
all one-form $a$. Since, for all p-forms $\psi$, $\phi$, we have
$\psi \wedge \ast \varphi = \langle \psi , \phi \rangle {\sl
Vol}$, it follows part i) of next lemma.
\begin{lemma} \label{astperpl}
{\it An almost quaternion-Hermitian manifold satisfies:
\begin{enumerate}
\item[{\rm i)}] For all five-form $\varphi$, $\varphi \in \left(
{\sl E} ({\sl H} + {\sl S}^3 {\sl H} + {\sl S}^5 {\sl H})
\right)^{\perp}$ if and only if, for $A,B =I,J,K$, $ \ast \varphi
\wedge \omega_{A} \wedge \omega_{B}
 = 0$.
\item[{\rm ii)}] For $A=I,J,K$, we have $2 \langle A \cdot
\lrcorner {\sl d^*} \omega_A, \omega_A \rangle = \ast ( \ast {\sl
d} \Omega \wedge \omega_{A} \wedge \omega_{A})$.
\end{enumerate}
 }
\end{lemma}
 \Dt It remains to prove ii). For all three one-forms $\zeta_I$, $\zeta_J$, $\zeta_K$ and $A =
I,J,K$, by direct computation we obtain \vspace{-6mm}

\begin{equation} \label{astff}
\qquad \ast \left(  \mbox{$\sum_{B,C=I,J,K}$} \ast \left( {\sl
i}_{B} ( C\zeta_{C} \wedge \omega_{C} ) \wedge \omega_{B} \right)
\wedge \omega_{A} \wedge \omega_{A} \right) = - 4 (n-1)(2n+1)
\zeta_{A}.
\end{equation}
Furthermore, we have  ${\sl d} \Omega = ({\sl d}\Omega)_1 + ({\sl
d} \Omega)_2$, where $({\sl d}\Omega)_1 \in (\Lambda_0^3 {\sl E} +
{\sl K}) ({\sl H} + {\sl S}^3 {\sl H})$ and $({\sl d} \Omega)_2
\in E ({\sl H} + {\sl S}^3 {\sl H})$. By part i), we have  $\ast
({\sl d} \Omega)_1 \wedge \omega_{A} \wedge \omega_{A}=0$. Now,
from Table 2,  taking \ref{xialpha}, $k_1=n-1$ and $k_2=2n+1$ into
account, we deduce that
\begin{equation} \label{es3heta}
 (  {\sl d} \Omega)_2 = \sum_{A,B=I,J,K} {\sl i}_{A} \left(
B\eta_{B} \wedge \omega_{B}
 \right)\wedge \omega_A,
\end{equation}
where the one-forms $\eta_{A}$ are given by $ \eta_{A} = -
\frac{1}{2k_1k_2} \langle A \cdot \lrcorner {\sl d} \omega_A,
\omega_A \rangle$. \vspace{1mm}

Hence,  taking \ref{astff} and \ref{es3heta}  into account, it
follows
$$
\ast \left( \ast {\sl d} \Omega \wedge \omega_{A} \wedge
\omega_{A} \right) = \ast \left( \ast ({\sl d} \Omega)_2 \wedge
\omega_{A} \wedge \omega_{A} \right) = 2 \langle A \cdot \lrcorner
{\sl d} \omega_A, \omega_A \rangle. \qedhere
$$

Next corollary is an immediate consequence of \ref{xi},
\ref{xialpha} and Lemma \ref{astperpl}.

\begin{proposition}
{\it \begin{enumerate} \item[{\rm i)}]$\pi_{{\sl E} {\sl H} }
(\nabla \Omega) = 0$ if and only if $ \ast {\sl d} \Omega \wedge
\Omega=0$.
 \item[{\rm ii)}]$\pi_{{\sl E}  {\sl S}^3 {\sl H}}
(\nabla \Omega) = 0$ if and only if $\ast {\sl d} \Omega \wedge
\omega_{I} \wedge \omega_{I}= \ast {\sl d} \Omega \wedge
\omega_{J} \wedge \omega_{J} =\ast {\sl d} \Omega \wedge
\omega_{K} \wedge \omega_{K}$.
 \item[{\rm iii)}]  $\pi_{{\sl E} ({\sl H} + {\sl S}^3 {\sl H})}
(\nabla \Omega) = 0$ if and only if $\ast {\sl d} \Omega \wedge
\omega_{A} \wedge \omega_{A}=0$, for $A=I,J,K$. \cqd
\end{enumerate}}
\end{proposition}
\vspace{2mm}

In the following lines we derive an expression for ${\sl d}^*
\Omega$ which is useful to handle examples. We recall the
following identity, given in \cite{Gray},
\begin{equation} \label{grayid}
 2 \nabla \omega_I =  d \omega_I
- I_{(2)} I_{(3)}  d \omega_I - I_{(2)} N_I,
\end{equation}
 where $N_I(X,Y,Z) = \langle X , N_I (Y, Z) \rangle$ and the (1,2)-tensor
$N_I$ is the Nijenhuis tensor for $I$ defined by $N_I(X,Y) =[X,Y]
+ I[IX,Y] + I[X,IY] - [IX,IY]$, for all vector fields $X$, $Y$,
$Z$. Moreover, we have the following fact, noted in
\cite{Alekseevski-Marchiafava},
\begin{equation} \label{nitra}
N_I(e_i,e_i,\cdot)=0.
\end{equation}
Now,  identities \ref{grayid} and \ref{nitra} imply
\begin{equation} \label{Leedd}
 I d^* \omega_I = - \langle \cdot \lrcorner d \omega_I, \omega_I
\rangle.
\end{equation}
On the other hand, from $ d^* \Omega = - \nabla_{e_i} \Omega (e_i,
\cdot, \cdot, \cdot) $ directly follows  the  expression  $ d^*
\Omega = 2 \sum_{A=I,J,K} \left( d^* \omega_A \wedge \omega_A - A
d \omega_A \right)$. Thus, taking \ref{Leedd} into account, we
obtain
\begin{equation} \label{codeOme}
{\sl d}^* \Omega = - 2 \sum_{A=I,J,K} \langle A \cdot \lrcorner d
\omega_A, \omega_A \rangle \wedge \omega_A - 2\sum_{A=I,J,K} A d
\omega_A.
\end{equation}

Finally, we briefly point out  a relation of examples already
indicated in some references and others studied by ourselves.
Because it would take up a lot of space in the present exposition,
we reserve a more detailed presentation of these latter examples
to future paper. If $\nabla \Omega=0$, the manifold is said to be
quaternionic K{\"a}hler (q.K.) and its metric $\langle,\rangle$ is
Einstein. The model example of such manifolds is the quaternionic
projective space ${\cal H}P(n)$. For more q.K. examples, there is
a relative extensive bibliography of them (for instance,
\cite{S1}). If $\nabla \Omega \in {\sl E}{\sl H}$, then the
manifold is said to be locally conformal quaternionic K{\"a}hler
(l.c.q.K.). The manifolds $S^{4n+3} \times S^1$ are locally
conformal hyperK{\"a}hler which is an special case of l.c.q.K. More
examples of these manifolds are given by L. Ornea and P. Piccinni
in \cite{OP}.
 The theory of quaternionic manifolds have been
independently developed by S. Salamon in \cite{S3} and by L.
Berard Bergery and T. Ochoiai in \cite{BO}. An almost
quaternion-Hermitian  (a.q.H.) manifold  is quaternionic if and
only if $\nabla \Omega \in (\Lambda_0^3 {\sl E} + {\sl K} + {\sl
E}) {\sl H}$ (see \cite{S3}). Compact examples of quaternionic
manifolds are given by D. Joyce in \cite{J}.

An example of eight-dimensional a.q.H. manifold  such that ${\sl
d} \Omega = 0$ and \linebreak $\nabla \Omega \neq 0$ (${\sl K}
{\sl S}^3 {\sl H}$ type) is given by S. Salamon in \cite{S4}.
Likewise one can check that the manifolds $N_1$, $N_2$ and $N_3$
given by I. Dotti and A. Fino in \cite{DF} are  a.q.H. manifolds
of type ${\sl K}{\sl H}$. The quaternionic Heisenberg group
studied  by L. Cordero, M. Fern{\'a}ndez and M. de Le{\'o}n in \cite{CFL},
is an a.q.H. manifold of type $(\Lambda_0^3 {\sl E} + {\sl K})
{\sl H}$. Quaternion K{\"a}hler manifolds with torsion, introduced by
S. Ivanov in \cite{Ivanov:QKT}, can be identified with the class
$({\sl K} + {\sl E}){\sl H}$. Special cases of these manifolds are
hyperK{\"a}hler manifolds with torsion (HKT). The manifold $({\sl S}^1
\times {\sl S}^3)^3$ is HKT and hence $({\sl K}   + {\sl E}) {\sl
H}$. Because the torsion one-form is closed, doing local conformal
changes of metric one can obtain open submanifolds of type ${\sl
K} {\sl H}$. On $S^3 \times T^9$ and $\left( {\sl S}^3 \right)^4$,
one can find  a.q.H. structures of type $\Lambda_0^3 {\sl E} ({\sl
S}^3 {\sl H} + {\sl H})$. On $T^3 \times M^3$, with $M$ a
three-dimensional Lie group, either nilpotent or solvable, one can
define a.q.H. structures of type $({\sl K}+{\sl E})({\sl S}^3 {\sl
H} + {\sl H})$, $(\Lambda_0^3 {\sl E} + {\sl K}) ( {\sl S}^3 {\sl
H} + {\sl H})$, $({\sl K} + {\sl E}){\sl S}^3 {\sl H}$ and
$(\Lambda_0^3 {\sl E} + {\sl K}){\sl H} + {\sl K} {\sl S}^3 {\sl
H}$. In the case of $T^3 \times M^3$ with a.q.H. of type $({\sl
K}+{\sl E})({\sl S}^3 {\sl H} + {\sl H})$, the one-form $\xi$ is
closed. Therefore, doing local conformal changes of metric, one
can obtain open submanifolds with a.q.H. structure of type ${\sl
K}({\sl S}^3 {\sl H} + {\sl H}) + {\sl E} {\sl S}^3 {\sl H}$.

\begin{observation}
 The results here contained  are also valid for almost quaternion
pseudo-Hermitian manifolds as it happens with Swann's results in
\cite{Sw1}. In such a case, we have an  ${\sl Sp}(p,q){\sl
Sp}(1)$-structure and $(4p,4q)$ is the signature of the metric
$\langle,\rangle$. Since $4q$ is even the expressions involving
the Hodge $\ast$-operator still are valid. In the present text, we
would only need to do slight modifications in those expressions
involving an orthonormal basis $\{ e_1, \dots, e_{4n} \}$ for
vectors, where we would have to write $\epsilon_r = \langle e_r ,
e_ r \rangle$. For instance, the complex symplectic form
$\omega_E$ and the extension of~$\langle\cdot ,\cdot \rangle$ to
$(0,s)$-tensors would be respectively given by
\begin{displaymath}
 \omega_{E} = \epsilon_i e^{\omega}_i \wedge
\widetilde{e}^{\omega}_i, \quad   \langle a, b \rangle =
\mbox{$\frac{1}{s!}$} \, a (
  e_{i_1}, \dots , e_{i_s})  b (e_{i_1}, \dots , e_{i_s}) \epsilon_{i_1} \dots
  \epsilon_{i_s}.
\end{displaymath}
\end{observation}

\newpage
{\scriptsize
\begin{center}
{\normalsize Table 2: Classes of almost quaternion-Hermitian
manifolds of dimension $\geq 12$}
  \label{tados}
\vspace{1mm}

\begin{tabular}{|p{1.9cm}|p{6.2cm}|p{5.3cm}|}
\hline
 & & \\[-3mm]
{\scriptsize  $ {\cal QK}$}& {\scriptsize  $\nabla \Omega  = 0$}
& {\scriptsize $\mbox{\sl d} \Omega = 0$}  \\[1mm]
\hline
 &&\\[-3mm]
{\scriptsize  $\Lambda_0^3 {\sl E} {\sl H}$}& {\scriptsize
$\displaystyle {\cal L}( \nabla \Omega ) = 4 \nabla \Omega $ and
$\mbox{\sl d}^* \Omega = 0$} & {\scriptsize $\mbox{\sl L}(
\mbox{\sl d} \Omega) = 6 d \Omega $ and $d^* \Omega = 0$ }
\\[1mm]
\hline
 &&\\[-3mm]
{\scriptsize  ${\sl K} {\sl H}$}&  {\scriptsize $ \displaystyle
\nabla_{\cdot} \Omega =
 \frac{1}{6} \sum_{A=I,J,K} \mbox{\sl i}_A (\cdot \lrcorner d^*
\Omega)  \wedge \omega_A $ \newline and  $\xi_I = \xi_J =\xi_K $
}& {\scriptsize $ \displaystyle  \mbox{\sl d} \Omega = \frac{1}{3}
\sum_{A=I,J,K} \mbox{\sl i}_A (\mbox{\sl d}^* \Omega) \wedge
\omega_A$ \newline and  $\xi_I = \xi_J =\xi_K $ }
\\[1mm]
\hline
 &&\\[-3mm]
{\scriptsize  ${\sl E} {\sl H}$}& {\scriptsize $ \displaystyle
\nabla_{\cdot} \Omega =  - \frac{1}{4k_1} \{  \cdot \wedge ( \xi
\lrcorner \Omega) - \xi \wedge ( \cdot \lrcorner  \Omega) \} $ } &
{\scriptsize $\displaystyle \mbox{\sl d} \Omega = - \frac{1}{k_1}
\xi \wedge \Omega$}
\\[1mm]
\hline
 & &\\[-3mm]
{\scriptsize  $\Lambda_0^3 {\sl E} {\sl S}^3 {\sl H}$}&
{\scriptsize
 $ \displaystyle \nabla_{\cdot}
\Omega = -  \frac{1}{6} \sum_{A=I,J,K} \mbox{\sl i}_A (\cdot
\lrcorner d^* \Omega)  \wedge \omega_A$ and $ \xi = 0$ } &
{\scriptsize $\displaystyle \mbox{\sl d} \Omega = - \frac{1}{3}
\sum_{A=I,J,K} \mbox{\sl i}_A (\mbox{\sl d}^* \Omega) \wedge
\omega_A$ \newline and $ \xi = 0$ }
\\[1mm]
\hline
 &&\\[-3mm]
{\scriptsize  ${\sl K} {\sl S}^3 {\sl H}$}&  {\scriptsize $
\displaystyle {\cal L}(\nabla \Omega)  = - 2 \nabla \Omega$ and
$\mbox{\sl d}^* \Omega = 0$  } & {\scriptsize $ \mbox{\sl L}
(\mbox{\sl d} \Omega ) = 0$ and $\mbox{\sl d}^* \Omega = 0$ }
\\[1mm]
\hline
 & &\\[-3mm]
{\scriptsize  ${\sl E} {\sl S}^3 {\sl H}$}& {\scriptsize
 $
\displaystyle \nabla_{\cdot} \Omega =  \frac{1}{k_2}
\sum_{A,B=I,J,K} \mbox{\sl i}_A \left( \cdot \lrcorner
\left(B\xi_{B} \wedge \omega_{B} \right)\right) \wedge \omega_{A}
$ and $\xi = 0$} & {\scriptsize $ \displaystyle \mbox{\sl d}
\Omega = \frac{2}{k_2} \sum_{A,B=I,J,K} \mbox{\sl i}_A \left(
B\xi_B \wedge \omega_B \right) \wedge \omega_A$
\newline and $\xi=0$ }
\\[1mm]
\hline
 & &\\[-3mm]
{\scriptsize  $ ( \Lambda_0^3 {\sl E}  +  {\sl K}) {\sl H} $}&
{\scriptsize $\displaystyle {\cal L} (\nabla \Omega ) = 4 \nabla
\Omega$ and  $\xi =0$} & {\scriptsize $\mbox{\sl L}(\mbox{\sl d}
\Omega) = 6 \mbox{\sl d} \Omega$ and $\xi = 0$ }
 \\[1mm]
\hline
 & & \\[-3mm]
{\scriptsize  $ (\Lambda_0^3 {\sl E}  +  {\sl E}) {\sl H} $}&
{\scriptsize ${\cal L} (\nabla \Omega ) = 4 \nabla \Omega$ and
$\mbox{\sl d}^* \Omega = \displaystyle \xi \lrcorner \Omega$} &
{\scriptsize$\mbox{\sl L}(\mbox{\sl d} \Omega) = 6 \mbox{\sl d}
\Omega$  and $\mbox{\sl d}^* \Omega = \displaystyle \xi \lrcorner
\Omega$}
 \\[1mm]
\hline
 & &\\[-3mm]
{\scriptsize  $ \Lambda_0^3 {\sl E} ( {\sl H} +  {\sl S}^3 {\sl H}
)$}& {\scriptsize ${\cal L} (\nabla_{\cdot} \Omega ) = 4 \nabla_{\cdot} \Omega
+ \displaystyle \sum_{A=I,J,K} \mbox{\sl i}_A (\cdot \lrcorner d^*
\Omega)  \wedge \omega_A $ \newline and $\xi =0$} & {\scriptsize
$\displaystyle \mbox{\sl L} (\mbox{\sl d} \Omega) = 6 \mbox{\sl d}
\Omega + 2 \sum_{A=I,J,K} \mbox{\sl i}_A \left( \mbox{\sl d}^*
\Omega \right) \wedge \omega_A $ and $\xi=0$}
\\[1mm]
\hline
 &&\\[-3mm]
{\scriptsize  $   \Lambda_0^3 {\sl E}  {\sl H} + {\sl K} {\sl S}^3
{\sl H}$}&  {\scriptsize $\mbox{\sl d}^* \Omega = 0$  or
$\Omega^{n-2} \wedge \mbox {\sl d} \Omega = 0$ } & {\scriptsize
Idem}
 \\[1mm]
\hline
 & & \\[-3mm]
{\scriptsize  $ \Lambda_0^3 {\sl E} {\sl H} +  {\sl E}{\sl S}^3
{\sl H} $}& {\scriptsize $\displaystyle {\cal L}( \nabla_{\cdot} \Omega
) = 4 \nabla_{\cdot} \Omega \newline  \mbox{ } \hspace{.6cm} -
\frac{6}{k_2} \sum_{A,B=I,J,K} \mbox{\sl i}_A \left(\cdot
\lrcorner  \left ( B \xi_{B} \wedge \omega_{B} \right)\right)
\wedge \omega_{A}$
\newline and  $\mbox{\sl d}^* \Omega = \displaystyle - 2
\sum_{A=I,J,K} A \xi_{A}  \wedge \omega_{A}$}
 & {\scriptsize $ \displaystyle \mbox{\sl L} ( \mbox{\sl d} \Omega ) = 6 \mbox{\sl d} \Omega \newline \mbox{ } \hspace{5mm} -
 \frac{12}{k_2}\sum_{A,B=I,J,K} \mbox{\sl i}_A \left( B\xi_B \wedge \Omega_B
 \right) \wedge \omega_A$ \newline and $\mbox{\sl d}^* \Omega =
\displaystyle - 2 \sum_{A=I,J,K} A \xi_{A}  \wedge \omega_{A}$}
 \\[1mm]
\hline
 && \\[-3mm]
{\scriptsize  $({\sl K}  + {\sl E}){\sl H}$}& {\scriptsize   $
\displaystyle \nabla_{\cdot} \Omega =  \frac{1}{6} \sum_{A=I,J,K}
\mbox{\sl i}_A (\cdot \lrcorner d^* \Omega)   \wedge \omega_{A}
\newline \mbox{ } \hspace{.6cm} - \frac{k_2}{12k_1} \{ \cdot \wedge ( \xi \lrcorner
\Omega) - \xi \wedge ( \cdot \lrcorner  \Omega) \}$
\newline and $ \xi_I = \xi_J = \xi_K $   } & {\scriptsize
$\displaystyle \mbox{\sl d} \Omega = \frac{1}{3} \sum_{A=I,J,K}
\mbox{\sl i}_A ( \mbox{\sl d}^* \Omega ) \wedge \omega_A -
\frac{k_2}{3k_1} \xi \wedge \Omega $ and $\xi_I= \xi_J =
\xi_K$}
\\[1mm]
\hline
&&\\[-3mm]
{\scriptsize ${\sl K}{\sl H} + \Lambda_0^3 {\sl E} {\sl S}^3 {\sl
H} $}&
 {\scriptsize $ \displaystyle \nabla_{\cdot} \Omega =  \frac{1}{18} \sum_{A=I,J,K}  \mbox{\sl i}_{A} ( \cdot
\lrcorner \mbox{\sl L}(\mbox{\sl d}^* \Omega) )  \wedge \omega_{A}
$ } & {\scriptsize $ \displaystyle \mbox{\sl d} \Omega =
\frac{1}{9} \sum_{A=I,J,K} \mbox{\sl i}_A  ( \mbox{\sl
L}(\mbox{\sl d}^* \Omega) ) \wedge \omega_A$ }
\\[1mm]
\hline
& &\\[-2mm]
{\scriptsize $\mbox{\sl K} (\mbox{\sl H} + {\sl S}^3 {\sl H}) $}&
{\scriptsize $\displaystyle {\cal L} (\nabla_{\cdot}  \Omega ) = - 2
\nabla_{\cdot} \Omega + \displaystyle \sum_{A=I,J,K} \mbox {\sl i}_{A}
( \cdot \lrcorner \mbox{\sl d}^* \Omega )  \wedge \omega_{A}$
\newline
and $\xi_I = \xi_J = \xi_K $} & {\scriptsize $\mbox{\sl L}
(\mbox{\sl d} \Omega) = 2 \displaystyle \sum_{A=I,J,K} \mbox{\sl
i}_A(\mbox{\sl d}^* \Omega) \wedge \omega_A$ \newline and $\xi_I =
\xi_J = \xi_K$}
\\[1mm]
\hline
 & & \\[-3mm]
{\scriptsize ${\sl K} {\sl H} +  {\sl E}{\sl S}^3{\sl H}$}&
{\scriptsize $ \displaystyle \nabla_{\cdot} \Omega =  \frac{1}{6}
\sum_{A=I,J,K} \mbox{\sl i}_{A}  ( \cdot \lrcorner \mbox{\sl d}^*
\Omega )  \wedge \omega_{A} \newline \mbox{ } \hspace{.5cm} +
\frac{k_2 + 3}{3 k_2} \sum_{A,B=I,J,K}  \mbox{\sl i}_{A}  ( \cdot
\lrcorner (B\xi_{B} \wedge \omega_{B} )) \wedge \omega_{A}$} &
{\scriptsize $\displaystyle \mbox{\sl d} \Omega = \frac{1}{3}
\sum_{A=I,J,K} \mbox{\sl i}_A(\mbox{\sl d}^* \Omega) \wedge
\omega_A \newline + \frac{2(k_2+3)}{3k_2} \sum_{A,B=I,J,K}
\mbox{\sl i}_A (B \xi_B \wedge \omega_B) \wedge \omega_A $ }
\\
 \hline
\end{tabular}
\end{center}
}

\newpage
{\scriptsize
\begin{center}
\begin{tabular}{|p{2cm}|p{6.2cm}|p{5.3cm}|}
\hline
 & &\\[-3mm]
{\scriptsize $ {\sl E}{\sl H} + \Lambda_0^3 {\sl E} {\sl S}^3 {\sl
H}$}& {\scriptsize $ \displaystyle \nabla_{\cdot} \Omega = - \frac{1}{6}
\sum_{A=I,J,K} \mbox{\sl i}_{A} ( \cdot \lrcorner \mbox{\sl d}^*
\Omega )  \wedge \omega_{A} \newline \mbox{ } \hspace{.8cm} +
\frac{k_2-6}{12 k_1} \left\{  \cdot \wedge ( \xi \lrcorner \Omega)
- \xi \wedge ( \cdot \lrcorner  \Omega) \right\} $ }& {\scriptsize
$ \displaystyle \mbox{\sl d} \Omega = - \frac{1}{3} \sum_{A=I,J,K}
\mbox{\sl i}_{A} (  \mbox{\sl d}^* \Omega )  \wedge \omega_{A}
\newline \mbox{ } \hspace{.6cm} + \frac{k_2-6}{3 k_1} \xi \wedge
\Omega $ }
 \\[1.5mm]
\hline
 & & \\[-2mm]
{\scriptsize $   {\sl E}  {\sl H} +  {\sl K} {\sl S}^3 {\sl H}
$}&{\scriptsize $\displaystyle {\cal L} (\nabla_{\cdot} \Omega ) =
- 2 \nabla_{\cdot} \Omega - \frac{3}{2 k_1} \left\{ \cdot \wedge (
\xi \lrcorner \Omega) - \xi  \wedge ( \cdot \lrcorner \Omega)
\right\}$ \newline and $\mbox{\sl d}^*\Omega = \xi \lrcorner
\Omega$ } & {\scriptsize $\displaystyle \mbox{\sl L} ( \mbox{\sl
d} \Omega ) = - \frac{6}{ k_1} \xi \wedge  \Omega$ and $\mbox{\sl
d}^*\Omega = \xi \lrcorner \Omega$ }
 \\[1mm]
\hline
 & &\\[-3mm]
{\scriptsize $   {\sl E}{\sl H} + {\sl E} {\sl S}^3 {\sl H}$}&
{\scriptsize  $ \displaystyle \nabla_{\cdot} \Omega =
\frac{1}{k_2} \sum_{A,B=I,J,K}  {\sl i}_{A} ( \cdot \lrcorner (B
\xi_{B} \wedge \omega_{B} )) \wedge \omega_{A}  \newline \mbox{ }
\hspace{.6cm} - \frac{3}{4 k_1 k_2}
 \left\{  \cdot \wedge ( \xi \lrcorner  \Omega)
- \xi \wedge ( \cdot \lrcorner  \Omega) \right\} $ } &
{\scriptsize  $ \displaystyle \mbox{\sl d} \Omega = \frac{2}{k_2}
\sum_{A,B=I,J,K}  {\sl i}_{A} ( B \xi_{B} \wedge \omega_{B})
\wedge \omega_{A}  \newline \mbox{ } \hspace{.6cm} - \frac{3}{ k_1
k_2}  \xi \wedge  \Omega $ }  \\[1.5mm]
\hline && \\[-3mm]
 {\scriptsize $  ( \Lambda_0^3 {\sl E} + {\sl K} ) {\sl S}^3 {\sl H}  $}&
 {\scriptsize
 $\displaystyle
 {\cal L} (\nabla \Omega ) = - 2 \nabla \Omega $ and $\xi_I= \xi_J = \xi_K
 $} & {\scriptsize
 $\displaystyle
 \mbox{\sl L} (\mbox{\sl d} \Omega ) = 0 $ and $\xi_I= \xi_J = \xi_K
 $}
 \\[1mm]
\hline
 & &\\[-3mm]
{\scriptsize $(\Lambda_0^3 {\sl E} + {\sl E}){\sl S}^3{\sl H}$}&
{\scriptsize  $ \displaystyle \nabla_{\cdot} \Omega = -
\frac{1}{6} \sum_{A=I,J,K} \mbox{\sl i}_{A}  ( \cdot \lrcorner
\mbox{\sl d}^* \Omega )  \wedge \omega_{A} \newline \mbox{ }
\hspace{.6cm} - \frac{2 k_1}{3k_2} \sum_{A,B=I,J,K} \mbox{\sl
i}_{A} ( \cdot \lrcorner (B \xi_{B} \wedge \omega_{B} )) \wedge
\omega_{A}$  and   $\xi = 0$} & {\scriptsize  $ \displaystyle
\mbox{\sl d} \Omega = - \frac{1}{3} \sum_{A=I,J,K} \mbox{\sl
i}_{A} (  \mbox{\sl d}^* \Omega )  \wedge \omega_{A}
\newline \mbox{ } \hspace{.6cm} - \frac{4 k_1}{3k_2}
\sum_{A,B=I,J,K} \mbox{\sl i}_{A} ( B \xi_{B} \wedge
\omega_{B} ) \wedge \omega_{A}$  and   $\xi = 0$}  \\[1mm]
\hline
 & &\\[-3mm]
{\scriptsize $ ({\sl K} +  {\sl E}) {\sl S}^3 {\sl H} $}&
{\scriptsize $\displaystyle {\cal L} ( \nabla_{\cdot} \Omega ) = -
2 \nabla_{\cdot} \Omega $ \newline
 and $\mbox{\sl d}^* \Omega = \displaystyle - 2
\sum_{A=I,J,K} A \xi_{A}  \wedge \omega_{A}$ }
 & {\scriptsize $\displaystyle \mbox{\sl L} ( \mbox{\sl d} \Omega ) =
 0 $
 \newline
 and $\mbox{\sl d}^* \Omega = \displaystyle - 2
\sum_{A=I,J,K} A \xi_{A}  \wedge \omega_{A}$ } \\[1mm]
\hline
 &&\\[-3mm]
{\scriptsize $ (\Lambda_0^3 {\sl E} + {\sl K} + {\sl E}) {\sl H}
$}& {\scriptsize $\displaystyle  {\cal L} ( \nabla \Omega ) = 4
\nabla \Omega $ } &  {\scriptsize $  \mbox{\sl L} ( \mbox{\sl d}
\Omega ) = 6 \mbox{\sl d} \Omega $ }
 \\[1mm]
\hline
 & & \\[-3mm]
{\scriptsize $ (\Lambda_0^3 {\sl E} +
  {\sl K}) {\sl H} +  \Lambda_0^3 {\sl E} {\sl S}^3 {\sl H}
$}& {\scriptsize $\displaystyle {\cal L} ( \nabla_{\cdot} \Omega )
= 4 \nabla_{\cdot} \Omega + \frac{1}{2} \sum_{A=I,J,K} \mbox{\sl
i}_{A}  ( \cdot \lrcorner \mbox{\sl d}^* \Omega )  \wedge
\omega_{A} \newline \mbox{ } \hspace{.4cm} - \frac{1}{6}
\sum_{A=I,J,K} {\sl i}_{A} ( \cdot \lrcorner \mbox{\sl
L}(\mbox{\sl d}^*  \Omega) ) \wedge \omega_{A}$ and $\xi = 0$ } &
{\scriptsize $\displaystyle \mbox{\sl L} ( \mbox{\sl d} \Omega ) =
6 \mbox{\sl d} \Omega +  \sum_{A=I,J,K} \mbox{\sl i}_{A} (
\mbox{\sl d}^* \Omega ) \wedge \omega_{A}
\newline \mbox{ } \hspace{.2cm} - \frac{1}{3} \sum_{A=I,J,K} {\sl
i}_{A} (  \mbox{\sl L}(\mbox{\sl d}^*  \Omega) ) \wedge
\omega_{A}$ and $\xi = 0$ } \\[-3mm]
\hline
 & &\\[-3mm]
{\scriptsize $  (\Lambda_0^3 {\sl E} + {\sl K}) {\sl H} + {\sl K}
{\sl S}^3 {\sl H} $}& {\scriptsize  $\mbox{\sl L}(\mbox{\sl d}^*
\Omega) = 3 \mbox{\sl d}^* \Omega$ and $\xi =0$} & {\scriptsize
Idem}
 \\[1mm]
\hline
 &  & \\[-3mm]
{\scriptsize $ (\Lambda_0^3 {\sl E} + {\sl K}){\sl H} + {\sl
E}{\sl S}^3 {\sl H}  $}& {\scriptsize $\displaystyle {\cal L} (
\nabla_{\cdot} \Omega ) =  4 \nabla_{\cdot}  \Omega -
\frac{6}{k_2} \sum_{A,B=I,J,K} {\sl i}_{A} ( \cdot \lrcorner (B
\xi_{B}  \wedge \mbox{ } \hspace{1.3cm} \omega_{B} )) \wedge
\omega_{A}$} & {\scriptsize $\displaystyle \mbox{\sl L} (
\mbox{\sl d} \Omega ) =  6 \mbox{\sl d}  \Omega - \frac{12}{k_2}
\sum_{A,B=I,J,K} {\sl i}_{A} ( B \xi_{B}  \wedge
\mbox{ } \hspace{1.2cm} \omega_{B} ) \wedge \omega_{A}$} \\[1mm]
\hline
 & &\\[-3mm]
{\scriptsize $ \Lambda_0^3 {\sl E}(H + {\sl S}^3 {\sl H})
\newline + {\sl E}{\sl H} $}& {\scriptsize $\displaystyle {\cal L} ( \nabla_{\cdot} \Omega ) =  4
\nabla_{\cdot} \Omega + \sum_{A=I,J,K}  {\sl i}_{A}  ( \cdot
  \lrcorner \mbox{\sl d}^* \Omega) \wedge \omega_{A} \newline
\mbox{ } \hspace{1cm} - \cdot \wedge ( \xi \lrcorner  \Omega) +
\xi \wedge ( \cdot \lrcorner \Omega)$} & {\scriptsize
$\displaystyle \mbox{\sl L} ( \mbox{\sl d} \Omega ) =  6 \mbox{\sl
d} \Omega + 2 \sum_{A=I,J,K} {\sl i}_{A} (  \mbox{\sl d}^* \Omega)
\wedge \omega_{A}
\newline \mbox{ } \hspace{1cm} - 4 \xi \wedge  \Omega$}  \\
\hline
 & &\\[-3mm]
{\scriptsize $    (\Lambda_0^3 {\sl E} + {\sl E})
 {\sl H} + {\sl K}  {\sl S}^3 {\sl
H} $}& {\scriptsize $\mbox{\sl d}^* \Omega = \xi \lrcorner \Omega$
}
& {\scriptsize Idem} \\[1mm]
\hline
 & &\\[-3mm]
{\scriptsize $ (\Lambda_0^3 {\sl E} + {\sl E})  {\sl H} + {\sl E}
{\sl S}^3 {\sl H} $}&{\scriptsize $\displaystyle {\cal L} (
\nabla_{\cdot} \Omega ) =  4 \nabla_{\cdot} \Omega -  \frac{6}{k_2}
\sum_{A,B=I,J,K} {\sl i}_{A} ( \cdot \lrcorner (B \xi_{B}  \wedge
\newline \mbox{ } \hspace{1cm} \omega_{B}) ) \wedge \omega_{A} -
\frac{3}{k_2} \{ \cdot \wedge (\xi \lrcorner \Omega) - \xi \wedge
( \cdot \lrcorner \Omega) \}$
\newline and $ \mbox{\sl d}^* \Omega = -2 \sum_{A=I,J,K} A \xi_A \wedge
\omega_A$ } & {\scriptsize $\displaystyle \mbox{\sl L} ( \mbox{\sl
d} \Omega ) = 6 \mbox{\sl d} \Omega - \frac{12}{k_2}
\sum_{A,B=I,J,K} {\sl i}_{A} ( B \xi_{B}  \wedge \newline \mbox{ }
\hspace{1cm} \omega_{B} ) \wedge \omega_{A} - \frac{12}{k_2} \xi
\wedge \Omega$ \newline and $ \mbox{\sl d}^* \Omega = -2
\sum_{A=I,J,K} A \xi_A \wedge \omega_A$ }
 \\[1mm]
\hline
 & &\\[-3mm]
{\scriptsize $   \Lambda_0^3 {\sl E} (  {\sl H} + {\sl S}^3 {\sl
H})  \newline + {\sl K} {\sl S}^3 {\sl H} $ }& {\scriptsize
$\mbox{\sl L}(\mbox{\sl d}^* \Omega) = - 3 \mbox{\sl d}^* \Omega$
and $\xi_I = \xi_J = \xi_K $} & {\scriptsize Idem}
 \\[1mm]
\hline
\end{tabular}
\end{center}
}

\newpage
{\scriptsize
\begin{center}
\begin{tabular}{|p{2cm}|p{6.2cm}|p{5.3cm}|}
\hline
 & &\\[-3mm]
{\scriptsize $ \Lambda_0^3{\sl E}({\sl H} + {\sl S}^3{\sl H} )
\newline + {\sl E}{\sl S}^3{\sl H}$}& {\scriptsize  $\displaystyle {\cal
L} ( \nabla_{\cdot} \Omega ) =  4 \nabla_{\cdot} \Omega +
\sum_{A=I,J,K}  \mbox{\sl i}_{A}  ( \cdot \lrcorner \mbox{\sl d}^*
\Omega )  \wedge \omega_{A} \newline + \frac{4 k_1 }{k_2}
\sum_{A,B=I,J,K} {\sl i}_{A}  ( \cdot \lrcorner (B\xi_{B} \wedge
\omega_{B} )) \wedge \omega_{A} $ \newline and  $\xi  = 0 $ } &
{\scriptsize $\displaystyle \mbox{\sl L} ( \mbox{\sl d} \Omega ) =
6 \mbox{\sl d} \Omega + 2 \sum_{A=I,J,K}  \mbox{\sl i}_{A}  (
\mbox{\sl d}^* \Omega )  \wedge \omega_{A} \newline + \frac{8 k_1
}{k_2} \sum_{A,B=I,J,K} {\sl i}_{A}  ( B\xi_{B} \wedge \omega_{B}
) \wedge \omega_{A} $ \newline and  $\xi  = 0 $ }
 \\[1mm]
\hline
 & &\\[-3mm]
{\scriptsize $      ( {\sl K} + {\sl E}) {\sl S}^3 {\sl H}
\newline
+ \Lambda_0^3 {\sl E}  {\sl H} $}& {\scriptsize
$\mbox{\sl d}^* \Omega = \displaystyle - 2 \sum_{A=I,J,K} A
\xi_{A} \wedge \omega_{A}  $
 and $\xi = 0$}& {\scriptsize Idem} \\[1mm]
\hline
 & & \\[-3mm]
{\scriptsize $    ({\sl K} + {\sl E}) {\sl H} \newline
 + \Lambda_0^3 {\sl  E}  {\sl S}^3 {\sl H}   $}& {\scriptsize $
\displaystyle \nabla_{\cdot} \Omega =  \frac{1}{18} \sum_{A=I,J,K}
\mbox{\sl i}_{A} ( \cdot \lrcorner \mbox{\sl L}(\mbox{\sl d}^*
\Omega) ) \wedge \omega_{A}    \newline \mbox{ } \hspace{.6cm} -
\frac{k_2}{12k_1} \{ \cdot \wedge (\xi \lrcorner  \Omega) - \xi
\wedge ( \cdot \lrcorner  \Omega) \}$} & {\scriptsize $
\displaystyle \mbox{\sl d} \Omega =  \frac{1}{9} \sum_{A=I,J,K}
\mbox{\sl i}_{A} (  \mbox{\sl L}(\mbox{\sl d}^* \Omega) ) \wedge
\omega_{A}    \newline \mbox{ } \hspace{.6cm} - \frac{k_2}{3k_1}
\xi \wedge \Omega $}\\[1mm]
\hline
 & &\\[-3mm]
{\scriptsize $  {\sl K} ({\sl H} + {\sl S}^3 {\sl H}) \newline +
{\sl E} {\sl H} $}&{\scriptsize $\displaystyle {\cal L}
(\nabla_{\cdot} \Omega ) = -2 \nabla_{\cdot} \Omega +
\sum_{A=I,J,K} {\sl i}_{A} ( \cdot \lrcorner \mbox{\sl d}^* \Omega
)  \wedge \omega_{A}
\newline \mbox{ } \hspace{1cm} - \frac{k_2 }{2k_1} \{ \cdot \wedge (\xi \lrcorner  \Omega)
- \xi \wedge ( \cdot \lrcorner  \Omega) \} $}&{\scriptsize
$\displaystyle \mbox{\sl L} (\mbox{\sl d} \Omega ) =   2
\sum_{A=I,J,K} {\sl i}_{A} (  \mbox{\sl d}^* \Omega ) \wedge
\omega_{A}
\newline \mbox{ } \hspace{1cm} - \frac{2k_2 }{k_1} \xi \wedge \Omega$}\\[1mm]
\hline
 & &\\[-3mm]
{\scriptsize ${\sl E}( {\sl H} + {\sl S}^3 {\sl H}) + {\sl K} {\sl
H} $}& {\scriptsize $\displaystyle \nabla_{\cdot} \Omega = \frac{1}{6}
\sum_{A=I,J,K}  {\sl i}_{A} ( \cdot \lrcorner \mbox{\sl d}^*
\Omega )  \wedge \omega_{A} \newline \mbox{ } \hspace{.6cm} +
\frac{k_2+3}{3k_2} \sum_{A,B=I,J,K} {\sl i}_{A} ( \cdot \lrcorner
(B\xi_{B} \wedge \omega_{B} )) \wedge \omega_{A}
\newline \mbox{ } \hspace{.6cm} - \frac{3}{4k_1k_2} \{\cdot \wedge
(\xi \lrcorner  \Omega) - \xi \wedge ( \cdot \lrcorner  \Omega)
\}$}& {\scriptsize $\displaystyle \mbox{\sl d} \Omega =
\frac{1}{3} \sum_{A=I,J,K} {\sl i}_{A} (  \mbox{\sl d}^* \Omega )
\wedge \omega_{A}
\newline  + \frac{2k_2+6}{3k_2}
\sum_{A,B=I,J,K} {\sl i}_{A} ( B\xi_{B} \wedge \omega_{B} ) \wedge
\omega_{A}
\newline   - \frac{3}{k_1k_2} \xi  \wedge
  \Omega$}\\[1mm]
\hline
 & & \\[-3mm]
{\scriptsize $  {\sl K} ({\sl H} + {\sl S}^3 {\sl H}) \newline +
\Lambda_0^3 {\sl E} {\sl S}^3 {\sl H} $}&{\scriptsize
$\displaystyle {\cal L} ( \nabla_{\cdot} \Omega ) = -2
\nabla_{\cdot}  \Omega + \frac{1}{2} \sum_{A=I,J,K} {\sl i}_{A} (
\cdot  \lrcorner {\sl d}^* \Omega )  \wedge \omega_{A} \newline
\mbox{ } \hspace{1cm}  + \frac{1}{6} \sum_{A=I,J,K} {\sl i}_{A}  (
\cdot \lrcorner {\sl L}({\sl d}^* \Omega) ) \wedge \omega_{A}$
\newline and  $\xi_I   = \xi_J = \xi_K $}& {\scriptsize
$\displaystyle {\sl L} ( {\sl d} \Omega ) =   \sum_{A=I,J,K} {\sl
i}_{A} (  {\sl d}^* \Omega ) \wedge \omega_{A}
\newline \mbox{ } \hspace{1cm}  + \frac{1}{3} \sum_{A=I,J,K} {\sl
i}_{A}  (  {\sl L}({\sl d}^* \Omega) ) \wedge
\omega_{A}$ \newline and  $\xi_I   = \xi_J = \xi_K $}\\[1mm]
\hline
 & & \\[-3mm]
{\scriptsize $(\Lambda_0^3 {\sl E} + {\sl E}) {\sl S}^3 {\sl H}
\newline + {\sl K}  {\sl H}  $}&{\scriptsize $ \displaystyle \nabla_{\cdot} \Omega =
\frac{1}{18} \sum_{A=I,J,K} {\sl i}_{A} ( \cdot  \lrcorner {\sl
L}({\sl d}^* \Omega) ) \wedge \omega_{A} \newline \mbox{ }
\hspace{.8cm} - \frac{2k_1}{3k_2} \sum_{A,B=I,J,K} {\sl i}_{A}  (
\cdot \lrcorner (B\xi_{B} \wedge \omega_{B} )) \wedge
\omega_{A}$}&{\scriptsize $ \displaystyle {\sl d} \Omega =
\frac{1}{9} \sum_{A=I,J,K} {\sl i}_{A} ({\sl L}({\sl d}^* \Omega)
) \wedge \omega_{A}
\newline \mbox{ } \hspace{.4cm} - \frac{4k_1}{3k_2} \sum_{A,B=I,J,K} {\sl i}_{A}  (B\xi_{B} \wedge
\omega_{B} ) \wedge \omega_{A}$} \\[1mm]
\hline
 & & \\[-3mm]
{\scriptsize $  {\sl K} ({\sl H} + {\sl S}^3 {\sl H}) \newline +
{\sl E} {\sl S}^3  {\sl H} $}&{\scriptsize  $\displaystyle {\cal
L}  ( \nabla_{\cdot} \Omega ) = -2 \nabla_{\cdot} \Omega + \sum_{A=I,J,K}
{\sl i}_{A}  ( \cdot \lrcorner {\sl d}^* \Omega ) \wedge
\omega_{A} \newline \mbox{ } \hspace{1.2cm} +  2 \sum_{A,B=I,J,K}
{\sl i}_{A} ( \cdot \lrcorner (B\xi_{B}  \wedge \omega_{B} ))
\wedge \omega_{A} $}&{\scriptsize  $\displaystyle {\sl L}  ( {\sl
d} \Omega ) = 2 \sum_{A=I,J,K} {\sl i}_{A} ( {\sl d}^* \Omega )
\wedge \omega_{A}
\newline \mbox{ } \hspace{.6cm} +  4 \sum_{A,B=I,J,K} {\sl i}_{A}
( B\xi_{B}  \wedge
\omega_{B} ) \wedge \omega_{A} $} \\[1mm]
\hline
 & &\\[-3mm]
{\scriptsize $   (\Lambda_0^3 {\sl E} + {\sl K} ) {\sl S}^3 {\sl
H} \newline +  {\sl E} {\sl H}  $}&{\scriptsize  $\displaystyle
{\cal L} ( \nabla_{\cdot} \Omega ) = -2 \nabla_{\cdot} \Omega -\frac{3}{2k_1}
\{  \cdot \wedge (\xi \lrcorner \Omega) - \xi \wedge ( \cdot
\lrcorner  \Omega) \}$ \newline and  $\xi_I = \xi_J = \xi_K $}&
{\scriptsize  $\displaystyle {\sl L} ( {\sl d} \Omega ) =
-\frac{6}{k_1} \xi \wedge \Omega$ and  $\xi_I = \xi_J = \xi_K $}
\\[1mm]
\hline
 &&\\[-3mm]
{\scriptsize ${\sl E}({\sl H} +{\sl S}^3{\sl H}) \newline +
\Lambda_0^3 {\sl E} {\sl S}^3 {\sl H} $}&{\scriptsize $
\displaystyle \nabla_{\cdot} \Omega =  - \frac{1}{6} \sum_{A=I,J,K} {\sl
i}_{A}  ( \cdot \lrcorner {\sl d}^*  \Omega ) \wedge \omega_{A}
\newline \mbox{ } \hspace{1cm} - \frac{2k_1}{3k_2}
\sum_{A,B=I,J,K} {\sl i}_{A}  ( \cdot \lrcorner (B \xi_{B}  \wedge
\omega_{B} )) \wedge \omega_{A} \newline \mbox{ } \hspace{1cm}  -
\frac{3}{4k_1k_2} \{ \cdot \wedge (\xi \lrcorner \Omega) - \xi
\wedge ( \cdot \lrcorner  \Omega) \} $}&{\scriptsize $
\displaystyle {\sl d} \Omega =  - \frac{1}{3} \sum_{A=I,J,K} {\sl
i}_{A}  (  {\sl d}^*  \Omega ) \wedge \omega_{A}
\newline \mbox{ } \hspace{.5cm} - \frac{4k_1}{3k_2}
\sum_{A,B=I,J,K} {\sl i}_{A}  ( B \xi_{B}  \wedge \omega_{B} )
\wedge \omega_{A} \newline \mbox{ } \hspace{.5cm}  -
\frac{3}{k_1k_2} \xi \wedge \Omega  $}  \\[1mm]
\hline
\end{tabular}
\end{center}
}

\newpage
{\scriptsize
\begin{center}
\begin{tabular}{|p{2cm}|p{6.2cm}|p{5.3cm}|}
\hline
 & &\\[-3mm]
{\scriptsize $   {\sl E}  ( {\sl H} +   {\sl S}^3 {\sl H})
\newline  + {\sl K}{\sl S}^3 {\sl H} $}&{\scriptsize $\displaystyle
{\cal L} ( \nabla_{\cdot}  \Omega ) = - 2 \nabla_{\cdot} \Omega -
\frac{3}{2k_1} \{  \cdot \wedge (\xi  \lrcorner \Omega) - \xi
\wedge ( \cdot \lrcorner  \Omega) \}$ and ${\sl d}^* \Omega =
\displaystyle - 2 \sum_{A=I,J,K} A\xi_{A}  \wedge
\omega_{A}$}&{\scriptsize $\displaystyle {\sl L} ( {\sl d} \Omega
) =  - \frac{6}{2k_1} \xi \wedge \Omega$ \newline  and ${\sl d}^*
\Omega = \displaystyle - 2 \sum_{A=I,J,K}
A\xi_{A}  \wedge \omega_{A}$}  \\[1mm]
\hline
 &&\\[-3mm]
{\scriptsize $(\Lambda_0^3 {\sl E} + {\sl K}+ {\sl E})
\newline {\sl S}^3 {\sl H} $}& {\scriptsize ${\cal L}( \nabla
\Omega ) = - 2 \nabla \Omega $}&{\scriptsize ${\sl L}( {\sl d}
\Omega ) = 0$}
  \\[1mm]
\hline
 & & \\[-3mm]
{\scriptsize $ \Lambda_0^3 {\sl E} ( {\sl H} + {\sl S}^3 {\sl H})
\newline + ({\sl K} + {\sl E}) {\sl H} $}&{\scriptsize $\displaystyle
 {\cal L} ( \nabla_{\cdot} \Omega) = 4
\nabla_{\cdot} \Omega + \frac{1}{2} \sum_{A=I,J,K} {\sl i}_{A}  (
\cdot \lrcorner {\sl d}^* \Omega )  \wedge \omega_{A}\newline
\mbox{ } \hspace{.8cm} - \frac{1}{6} \sum_{A=I,J,K} {\sl i}_{A}  (
\cdot \lrcorner {\sl L}({\sl d} \Omega) ) \wedge \omega_{A}
$}&{\scriptsize $\displaystyle
 {\sl L} ( {\sl d} \Omega) = 6
{\sl d} \Omega +  \sum_{A=I,J,K} {\sl i}_{A}  (  {\sl d}^* \Omega
)  \wedge \omega_{A}\newline \mbox{ } \hspace{.5cm} - \frac{1}{3}
\sum_{A=I,J,K} {\sl i}_{A}  (  {\sl L}({\sl d} \Omega) ) \wedge
\omega_{A}    $}  \\[1mm]
\hline
 & &\\[-3mm]
{\scriptsize $   {\sl K} ( {\sl H}  +  {\sl S}^3 {\sl H})\newline
+ (\Lambda_0^3 {\sl E} +  {\sl E}) {\sl H} $}&{\scriptsize ${\sl
L}({\sl d}^* \Omega) = 3 {\sl d}^* \Omega$} & {\scriptsize Idem}  \\[1mm]
\hline
 & & \\[-3mm]
{\scriptsize $  {\sl E} ( {\sl H} +  {\sl S}^3 {\sl H}) \newline +
(\Lambda_0^3 {\sl E} +  {\sl K}) {\sl H}$}&{\scriptsize
$\displaystyle {\cal L} ( \nabla_{\cdot} \Omega ) = 4
\nabla_{\cdot} \Omega - \frac{3}{k_2} \{ \cdot \wedge (\xi
\lrcorner \Omega) - \xi \wedge ( \cdot \lrcorner  \Omega) \}
\newline \mbox{ } \hspace{1cm}  - \frac{6}{k_2} \sum_{A,B=I,J,K} {\sl i}_{A} ( \cdot
\lrcorner (B\xi_{B} \wedge \omega_{B} )) \wedge \omega_{A}
$}&{\scriptsize $\displaystyle {\sl L} ( {\sl d} \Omega ) = 6 {\sl
d} \Omega - \frac{12}{k_2} \xi \wedge \Omega
\newline \mbox{ } \hspace{.6cm}  - \frac{12}{k_2} \sum_{A,B=I,J,K} {\sl i}_{A} (
B\xi_{B} \wedge \omega_{B} ) \wedge \omega_{A}
$} \\[1mm]
\hline
& &\\[-3mm]
{\scriptsize $ (\Lambda_0^3 {\sl E} +  {\sl K})\newline ( {\sl H}
+ {\sl S}^3  {\sl H})$}& {\scriptsize $\xi_I = \xi_J = \xi_K = 0$
\newline or $\ast {\sl d} \Omega \wedge \omega_{A} \wedge
\omega_{A} = 0$, for $A=I,J,K$}&{\scriptsize Idem}  \\[1mm]
\hline
& &\\[-3mm]
 {\scriptsize $ \Lambda_0^3 {\sl E} ( {\sl H} + {\sl S}^3 {\sl H})
 \newline + {\sl K} {\sl H} +   {\sl E} {\sl S}^3 {\sl H} $}&{\scriptsize  $\displaystyle
{\cal L} ( \nabla_{\cdot} \Omega ) = 4 \nabla_{\cdot} \Omega +
\frac{1}{2} \sum_{A=I,J,K} {\sl i}_{A}  ( \cdot \lrcorner {\sl d}
\Omega )  \wedge \omega_{A} \newline \mbox{ } \hspace{1cm}  -
\frac{1}{6} \sum_{A,B=I,J,K} {\sl i}_{A}  ( \cdot \lrcorner {\sl
L}({\sl d}^* \Omega) ) \wedge \omega_{A} \newline \mbox{ }
\hspace{1cm}  + \frac{4k_1}{k_2} \sum_{A,B=I,J,K} {\sl i}_{A} (
\cdot \lrcorner (B\xi_{B} \wedge \omega_{B} )) \wedge
\omega_{A}$}&{\scriptsize  $\displaystyle {\sl L} ( {\sl d} \Omega
) = 6 {\sl d} \Omega +  \sum_{A=I,J,K} {\sl i}_{A} (  {\sl d}
\Omega ) \wedge \omega_{A}
\newline \mbox{ } \hspace{.6cm}  - \frac{1}{3} \sum_{A,B=I,J,K}
{\sl i}_{A}  ( {\sl L}({\sl d}^* \Omega) ) \wedge \omega_{A}
\newline \mbox{ } \hspace{.6cm}  + \frac{8k_1}{k_2}
\sum_{A,B=I,J,K} {\sl i}_{A} ( B\xi_{B} \wedge
\omega_{B} ) \wedge \omega_{A}$}  \\[1mm]
\hline
 &&\\[-3mm]
{\scriptsize $ (\Lambda_0^3 {\sl E} + {\sl K}) {\sl H} \newline +
({\sl K} + {\sl E}) {\sl S}^3 {\sl H} $}&{\scriptsize
  ${\sl L}({\sl d}^* \Omega) =  3
{\sl d}^* \Omega + 12 \displaystyle \sum_{A=I,J,K} A \xi_{A}
\wedge \omega_{A}$}&{\scriptsize Idem}
 \\[1mm]
\hline
 & &\\[-3mm]
{\scriptsize $ (\Lambda_0^3 {\sl E} + {\sl E}) {\sl H} + \newline
 ( \Lambda_0^3 {\sl E}+  {\sl K}) {\sl S}^3 {\sl H} $}&{\scriptsize ${\sl L}({\sl d}^* \Omega) = - 3
{\sl d}^* \Omega + 6 \xi \lrcorner \Omega$ \newline   and $\xi_I =
\xi_J = \xi_K$}
 &{\scriptsize Idem}\\[1mm]
 \hline
 & & \\[-3mm]
{\scriptsize $ (\Lambda_0^3 {\sl E} +   {\sl E})   ({\sl H}
\newline + {\sl S}^3 {\sl H})$}&{\scriptsize $\displaystyle {\cal
L} ( \nabla_{\cdot}\Omega ) =  4 \nabla_{\cdot} \Omega  + \sum_{A=I,J,K}  {\sl
i}_{A}  ( \cdot \lrcorner {\sl d}^*  \Omega ) \wedge \omega_{A}
\newline \mbox{ } \hspace{1cm} - \frac{3}{k_2} \{ \cdot \wedge
(\xi \lrcorner \Omega) - \xi \wedge ( \cdot \lrcorner \Omega)
\}\newline \mbox{ } \hspace{1cm} + \frac{4k_1}{k_2}
\sum_{A,B=I,J,K}{\sl i}_{A} ( \cdot \lrcorner (B\xi_{B}  \wedge
\omega_{B} )) \wedge \omega_{A} $}& {\scriptsize $\displaystyle
{\sl L} ( {\sl d}\Omega ) =  6 {\sl d} \Omega +  2 \sum_{A=I,J,K}
{\sl i}_{A}  (  {\sl d}^* \Omega ) \wedge \omega_{A}   \newline
\mbox{ } \hspace{.6cm} - \frac{12}{k_2} \xi \wedge  \Omega
\newline \mbox{ } \hspace{.6cm} + \frac{8k_1}{k_2}
\sum_{A,B=I,J,K}{\sl i}_{A} ( B\xi_{B}  \wedge \omega_{B} ) \wedge \omega_{A} $} \\[1mm]
\hline
 & &\\[-3mm]
{\scriptsize $ (\Lambda_0^3 {\sl E} + {\sl E}) {\sl H}+
\newline ({\sl K} +  {\sl E}) {\sl
S}^3 {\sl H} $}&{\scriptsize
  ${\sl d}^* \Omega =  - 2 \displaystyle \sum_{A=I,J,K}
A\xi_{A}  \wedge \omega_{A}$}&{\scriptsize Idem}
 \\[1mm]
\hline
 &&\\[-3mm]
 $ \Lambda_0^3 {\sl E} ({\sl H} + {\sl
S}^3  {\sl H}) \newline +  ({\sl K} + {\sl E}) {\sl S}^3 {\sl H}
$&
  ${\sl L} ({\sl d}^* \Omega)  =  - 3 {\sl d}^* \Omega$& Idem   \\[1mm]
\hline
 & & \\[-3mm]
$ ({\sl K}+ {\sl E})  {\sl H} + \newline ( \Lambda_0^3  {\sl E} +
 {\sl K}) {\sl S}^3 {\sl H} $& $\displaystyle
{\cal L}( \nabla_{\cdot} \Omega ) = -2 \nabla_{\cdot} \Omega + \frac{1}{2}
\sum_{A=I,J,K} {\sl i}_{A}  ( \cdot \lrcorner {\sl d}^* \Omega )
\wedge \omega_{A}  \newline \mbox{ } \hspace{1cm} + \frac{1}{6}
\sum_{A=I,J,K} {\sl i}_{A}  ( \cdot \lrcorner {\sl L}({\sl d}^*
\Omega) ) \wedge \omega_{A}       \newline \mbox{ } \hspace{1cm} -
\frac{k_2}{2k_1} \left\{ \cdot \wedge (\xi \lrcorner  \Omega) -
\xi \wedge ( \cdot \lrcorner  \Omega) \right\}$\newline  and
$\xi_I = \xi_J = \xi_K $& $\displaystyle {\sl L}( {\sl d} \Omega )
=   \sum_{A=I,J,K} {\sl i}_{A} ({\sl d}^* \Omega ) \wedge
\omega_{A}
\newline \mbox{ } \hspace{1cm} + \frac{1}{3} \sum_{A=I,J,K} {\sl
i}_{A}  (  {\sl L}({\sl d}^* \Omega) ) \wedge \omega_{A}
\newline \mbox{ } \hspace{1cm} - \frac{2k_2}{k_1} \xi \wedge \Omega$
\newline and $\xi_I = \xi_J = \xi_K $
 \\[1mm]
\hline
\end{tabular}
\end{center}
}

\newpage

{\scriptsize
\begin{center}
\begin{tabular}{|p{2.3cm}|p{6.1cm}|p{5.1cm}|}
\hline
 & & \\[-3mm]
 $   ( {\sl K} +{\sl E}) {\sl H} + \newline (\Lambda_0^3
{\sl E}  +  {\sl E})  {\sl S}^3 {\sl H} $& $ \displaystyle
\nabla_{\cdot} \Omega =  \frac{1}{18} \sum_{A=I,J,K} {\sl i}_{A} (
\cdot \lrcorner {\sl L} ({\sl d}^*\Omega) ) \wedge \omega_{A}
 \newline \mbox{ } \hspace{.6cm} - \frac{4k_1^2 + k_2^2}{12k_1k_2} \{ \cdot
\wedge (\xi \lrcorner \Omega) - \xi \wedge ( \cdot \lrcorner
\Omega) \}
\newline \mbox{ } \hspace{.5cm} - \frac{2k_1}{3k_2} \sum_{A,B=I,J,K} {\sl
i}_{A} ( \cdot \lrcorner (B\xi_{B} \wedge \omega_{B} )) \wedge
\omega_{A}$& $ \displaystyle {\sl d} \Omega =  \frac{1}{9}
\sum_{A=I,J,K} {\sl i}_{A} (  {\sl L} ({\sl d}^*\Omega) ) \wedge
\omega_{A}
 \newline \mbox{ } \hspace{.4cm} - \frac{4k_1^2 + k_2^2}{3k_1k_2} \xi \wedge \Omega
\newline \mbox{ } \hspace{.4cm} - \frac{4k_1}{3k_2} \sum_{A,B=I,J,K} {\sl
i}_{A} ( B\xi_{B} \wedge \omega_{B} ) \wedge
\omega_{A}$ \\[1mm]
\hline
 & &\\[-3mm]
 $ ({\sl K} + {\sl E}) (  {\sl H}  + {\sl S}^3 {\sl
H})$& $\displaystyle {\cal L}( \nabla_{\cdot} \Omega ) = -2 \nabla_{\cdot}
\Omega  +
 \sum_{A=I,J,K}  {\sl
i}_{A}  ( \cdot \lrcorner {\sl d}^* \Omega )  \wedge \omega_{A}
\newline \mbox{ } \hspace{1cm} - \frac{3}{2k_1} \{ \cdot \wedge
(\xi \lrcorner \Omega) - \xi \wedge ( \cdot \lrcorner  \Omega) \}
\newline \mbox{ } \hspace{1cm} + 2 \sum_{ A,B=I,J,K} {\sl i}_{A}
( \cdot \lrcorner (B\xi_{B} \wedge \omega_{B} )  )  \wedge
\omega_{A}$&$\displaystyle {\sl L}({\sl d} \Omega ) = 2
 \sum_{A=I,J,K}  {\sl
i}_{A}  ( {\sl d}^* \Omega )  \wedge \omega_{A} \newline \mbox{ }
\hspace{.7cm} - \frac{6}{k_1} \xi \wedge \Omega
\newline \mbox{ } \hspace{.7cm} + 4 \sum_{ A,B=I,J,K} {\sl i}_{A}
( B\xi_{B}
\wedge \omega_{B}   )  \wedge \omega_{A}$  \\[1mm]
\hline
 && \\[-3mm]
 $ {\sl K} ( {\sl H} + {\sl S}^3  {\sl H}) + \newline
( \Lambda_0^3 {\sl E}  + {\sl E}) {\sl S}^3 {\sl H}$&
$\displaystyle {\cal L} ( \nabla_{\cdot} \Omega ) = -2 \nabla_{\cdot} \Omega +
\frac{1}{2} \sum_{A=I,J,K} {\sl i}_{A} ( \cdot \lrcorner {\sl d}^*
\Omega )  \wedge \omega_{A} \newline \mbox{ } \hspace{1cm} +
\frac{1}{6} \sum_{A=I,J,K} {\sl i}_{A}  ( \cdot \lrcorner {\sl
L}({\sl d}^* \Omega) ) \wedge \omega_{A}$&$\displaystyle {\sl L} (
{\sl d} \Omega ) =  \sum_{A=I,J,K} {\sl i}_{A} (  {\sl d}^* \Omega
) \wedge \omega_{A}
\newline \mbox{ } \hspace{1cm} + \frac{1}{3} \sum_{A=I,J,K} {\sl
i}_{A}  (  {\sl L}({\sl d}^* \Omega) ) \wedge
\omega_{A}$ \\[1mm]
\hline
 &\\[-3mm]
$ {\sl E} ( {\sl H} + {\sl S}^3 {\sl H})
 + \newline (\Lambda_0^3 {\sl E} +  {\sl K}) {\sl S}^3 {\sl
H}   $& $\displaystyle {\cal L} ( \nabla_{\cdot} \Omega ) =  -2
\nabla_{\cdot} \Omega - \frac{3}{2k_1} \{ \cdot \wedge ( \xi
\lrcorner \Omega) - \xi \wedge ( \cdot \lrcorner  \Omega) \}
$&$\displaystyle {\sl L} ( {\sl d} \Omega ) = - \frac{6}{2k_1} \xi
\wedge \Omega$  \\[1mm]
\hline
 &\\[-3mm]
$ (\Lambda_0^3 {\sl E}+{\sl K}) {\sl H} + \newline ( \Lambda_0^3
{\sl E} + {\sl K}) {\sl S}^3 {\sl H}  \newline + {\sl E}{\sl H} $&
 $\xi_I  = \xi_J = \xi_K $ or \newline $\ast {\sl d} \Omega \wedge \omega_I
 \wedge \omega_I = \ast {\sl d} \Omega \wedge \omega_J \wedge \omega_J
 =\ast {\sl d} \Omega \wedge \omega_K \wedge \omega_K$  & Idem  \\[1mm]
\hline
 & &\\[-3mm]
$ (\Lambda_0^3 {\sl E}+{\sl E}) {\sl H} + \newline ( \Lambda_0^3
{\sl E} + {\sl E}) {\sl S}^3 {\sl H}  \newline + {\sl K}{\sl H} $&
$\displaystyle {\cal L} ( \nabla_{\cdot} \Omega ) = 4 \nabla_{\cdot} \Omega +
\frac{1}{2} \sum_{A=I,J,K} {\sl i}_{A}  ( \cdot \lrcorner {\sl
d}^* \Omega )  \wedge \omega_{A}\newline \mbox{ } \hspace{.9cm} -
\frac{1}{6} \sum_{A=I,J,K} {\sl i}_{A}  ( \cdot \lrcorner {\sl
L}({\sl d}^* \Omega) ) \wedge \omega_A
 \newline \mbox{
} \hspace{.9cm} - \frac{2k_1}{k_2} \{ \cdot \wedge (\xi \lrcorner
\Omega) - \xi \wedge ( \cdot \lrcorner  \Omega) \}
\newline \mbox{ } \hspace{.9cm} +
\frac{4k_1}{k_2} \sum_{A,B=I,J,K}  {\sl i}_{A} ( \cdot \lrcorner
(B\xi_{B}  \wedge \omega_{B} )) \wedge \omega_{A} $&$\displaystyle
{\sl L} ( {\sl d} \Omega ) = 6 {\sl d} \Omega +  \sum_{A=I,J,K}
{\sl i}_{A}  ( {\sl d}^* \Omega ) \wedge \omega_{A}\newline \mbox{
} \hspace{.4cm} - \frac{1}{3} \sum_{A=I,J,K} {\sl i}_{A}  ( {\sl
L}({\sl d}^* \Omega) ) \wedge \omega_A
 \newline \mbox{
} \hspace{.4cm} - \frac{8k_1}{k_2} \xi \wedge \Omega
\newline \mbox{ } \hspace{.4cm} +
\frac{8k_1}{k_2} \sum_{A,B=I,J,K}  {\sl i}_{A} ( B\xi_{B}  \wedge \omega_{B} ) \wedge \omega_{A} $  \\[2mm]
\hline
 &&\\[-3mm]
 $({\sl K} + {\sl E})  ({\sl H} + {\sl S}^3 {\sl
H}) + \Lambda_0^3 {\sl EH}  $&
  $ {\sl L}({\sl d}^*  \Omega) = 3 {\sl d}^* \Omega + 6 \xi \lrcorner
\Omega + 12 \displaystyle \sum_{A=I,J,K} A\xi_{A}  \wedge
\omega_{A}  $ & Idem
 \\[1mm]
\hline
 & & \\[-3mm]
 $ (\Lambda_0^3 {\sl E} + {\sl K}) ( {\sl H} + {\sl S}^3 {\sl H}) +  {\sl E} {\sl S}^3 {\sl H}  $&
  $ \xi = 0 $  or $\ast {\sl d} \Omega \wedge \Omega = 0$
& Idem  \\[1mm]
\hline
 &&\\[-3mm]
$ (\Lambda_0^3 {\sl E} + {\sl E}) ({\sl H}
 + {\sl S}^3 {\sl H}) +  {\sl K} {\sl
S}^3 {\sl H}   $&
  $ {\sl L}({\sl d}^* \Omega ) = - 3 {\sl d}^* \Omega + 6 \xi \lrcorner \Omega $
& Idem \\[1mm]
\hline
 & &\\[-3mm]
$ ({\sl K} + {\sl E})( {\sl H} + {\sl S}^3 {\sl H})  + \Lambda_0^3
{\sl E} {\sl S}^3 {\sl H} \newline  $&$\displaystyle {\cal L} (
\nabla_{\cdot} \Omega ) =  -2 \nabla_{\cdot} \Omega  + \frac{1}{2}
\sum_{A=I,J,K} {\sl i}_{A}  ( \cdot \lrcorner {\sl d}^* \Omega )
\wedge \omega_{A}  \newline \mbox{ } \hspace{1cm} + \frac{1}{6}
\sum_{A=I,J,K} {\sl i}_{A} ( \cdot \lrcorner {\sl L}({\sl d}^*
\Omega) ) \wedge \omega_{A}  \newline \mbox{ } \hspace{1cm} -
\frac{k_2}{2k_1} \{ \cdot \wedge (\xi \lrcorner \Omega) - \xi
\wedge ( \cdot \lrcorner  \Omega) \} $&$\displaystyle {\sl L} (
{\sl d} \Omega ) =   \sum_{A=I,J,K} {\sl i}_{A} ( {\sl d}^* \Omega
) \wedge \omega_{A}
\newline \mbox{ } \hspace{1cm} + \frac{1}{3} \sum_{A=I,J,K} {\sl
i}_{A} (  {\sl L}({\sl d}^* \Omega) ) \wedge \omega_{A}
\newline \mbox{ } \hspace{1cm} - \frac{2k_2}{k_1} \xi \wedge \Omega $  \\[1mm]
\hline
 &\\[-3mm]
$(\Lambda_0^3  {\sl E} +  {\sl K} + {\sl E}) ({\sl H} \newline +
{\sl S}^3 {\sl H})   $&
 no relation & no relation  \\[1mm]
\hline
\end{tabular}
\end{center}
}
\vspace{4mm}

\newpage
{\scriptsize
\begin{center}
{\normalsize Table 3: Partial classification via ${\sl d} \Omega$,
for
 $4n = 8$.}

\vspace{1mm}

\begin{tabular}{|p{2.6cm}|p{9.6cm}|}
\hline
 &\\[-3mm]
 ${\sl K} {\sl S}^3 {\sl H}$& $
\displaystyle {\sl d} \Omega = 0$ \\[1mm]
\hline
 &\\[-3mm]
 ${\sl K} ( {\sl H} + {\sl S}^3 {\sl H})$& $\xi_I =
\xi_J =\xi_K = 0 $ or $  \ast {\sl d} \Omega \wedge \omega_{A}
\wedge \omega_{A}  =0$, for $A=I,J,K$
\\[1mm]
\hline
 &\\[-3mm]
 ${\sl E} {\sl H} + {\sl K} {\sl S}^3 {\sl H}$& $
\displaystyle {\sl d}
\Omega =  -    \xi \wedge   \Omega $    \\[1mm]
\hline
 &\\[-3mm]
$({\sl K} + {\sl E}) {\sl S}^3 {\sl H} $&
 $
\displaystyle {\sl d} \Omega =  \frac{2}{5} \sum_{A,B=I,J,K} {\sl
i}_{A}  (B \xi_{B}  \wedge \omega_{B} ) \wedge \omega_{A} $ and
$\xi = 0$, or ${\sl L}({\sl
d}\Omega)=0$    \\[1mm]
\hline
 &\\[-3mm]
 ${\sl K} ({\sl H} +  {\sl S}^3
{\sl H} ) + {\sl E}{\sl H}  $&
 $ \xi_I = \xi_J = \xi_K $ or ${\sl L}({\sl d}\Omega)= 6 {\sl d}\Omega$     \\[1mm]
\hline
 &\\[-3mm]
${\sl K} ( {\sl H} + {\sl S}^3 {\sl H}) + {\sl E}{\sl S}^3{\sl
H}$&   $ \xi = 0$ or $\ast {\sl d} \Omega \wedge \Omega =
0$     \\[1mm]
\hline
 &\\[-3mm]
  ${\sl E}({\sl H} + {\sl S}^3 {\sl H}) + {\sl K}
{\sl S}^3 {\sl H} $&
 $ \displaystyle {\sl d} \Omega =  \frac{2}{5} \sum_{A,B=I,J,K}   {\sl i}_{A}  (B \xi_{B}  \wedge
\omega_{B} ) \wedge \omega_{A} - \frac{3}{ 5}
  \xi \wedge   \Omega $   \\[1mm]
\hline
 &\\[-3mm]
  $    ({\sl K} + {\sl E}) ({\sl H} + {\sl S}^3 {\sl H})$& $ \displaystyle {\sl d} \Omega =
\frac{1}{3} \sum_{A=I,J,K} {\sl i}_{A} ({\sl d}^* \Omega) \wedge
\omega_{A}  + \frac{16}{15} \sum_{A,B=I,J,K} {\sl i}_{A} (B\xi_{B}
 \wedge \omega_{B} ) \wedge \omega_{A}
  - \frac{3}{5}   \xi \wedge \Omega  $  \\[2mm]
\hline
\end{tabular}
\end{center}

}

\baselineskip=5mm

\baselineskip=5mm {\footnotesize
\bibliographystyle{plain}
\addcontentsline{toc}{section}{References}
\bibliography{papers}
 }

\end{document}